\documentclass[mathpazo]{cicp}
\usepackage{graphicx}
\usepackage{amssymb}
\usepackage{epstopdf}
\usepackage{amssymb}
\usepackage[a4paper, total={6in, 9.1in}]{geometry}
\usepackage{multirow}
\usepackage{bm}
\usepackage{amsmath,amsfonts,amsthm,mathrsfs,amssymb,cite}
\usepackage{mathtools}
\usepackage{framed}
\usepackage{exscale}
\usepackage{relsize}
\usepackage{enumerate}
\usepackage{epstopdf}
\graphicspath{{pic/}} 
\usepackage{latexsym}
\usepackage[usenames]{color}

\usepackage{graphics}
\usepackage{framed}
\usepackage{enumerate}
\usepackage{hyperref}
\usepackage{xcolor}

\usepackage{amsmath}
\usepackage{tikz}

\hypersetup{
	colorlinks,
	linkcolor={blue!80!black},
	urlcolor={blue!80!black},
	citecolor={blue!80!black}
}
\usepackage[capitalize,noabbrev]{cleveref}

\numberwithin{equation}{section}

\title{An ALE numerical method with HLLC-2D solver for the two-phase flow ejecta transporting model}

\author[Jianqiao Zhang et.~al.]{Jianqiao Zhang\affil{1}, Wei Yan\affil{1}\corrauth, Xianggui Li\affil{2}}
\address{\affilnum{1}\ Department of Mathematics, Jilin University, Changchun, Jilin, China. \\       
           \affilnum{2}\ School of Applied Science, Beijing Information Science and Technology University, Bejing, China.}
\emails{{\tt jqzhang22@mails.jlu.edu.cn} (J.~Zhang), {\tt wyanmath@jlu.edu.cn} (W.~Yan),
         {\tt lixg@bistu.edu.cn} (X.~Li)}

\date{} 
\begin{document}
	
\begin{abstract}

This work presents an arbitrary Lagrangian Eulerian (ALE) method for the compressible two-phase flow ejecta transporting model with the HLLC-2D Riemann solver. We focus on researching the precise equation to describe the interactions between particle phase and flow phase. The calculation of the momentum and energy exchange across two phases is the key point during the procedure, which can be capable of maintaining the conservation of this system. For particles, tracking their trajectories within the mesh and elements is essential. Thereafter an ALE method instead of Lagrangian scheme is derived for the discretization of the equation to perform better with the complex motion of particles and flow. We apply the HLLC-2D Riemann solver to substitute the HLLC solver which relaxes the limitation for continuous fluxes along the edge. Meanwhile we propose a method for searching particles and provide a CFL-like condition based on this. Finally, we show some numerical tests to analysis the influence of particles on fluid and get a following effect between two phases. The model and the numerical method are validated through numerical tests to show its robustness and accuracy.

\medskip
	

\end{abstract}
\maketitle


\section{Introduction}\label{sec:intro}

The compressible multiphase flows in a gas-particle system is one of the most important problems in both nature and engineering systems. Due to the high-complexity and multi-scale natures of the gas-particle system, providing effective experimental methods and numerical simulations is essential in the research of gas-particle flow.

Generally in the early time, systematic research on two-phase flow began in the 1940s, but it was primarily focused on experimental observations and phenomenological descriptions. The term `Two-phase Flow' first appeared in the literature in \cite{3}. After that, people started to focus on this problem which led to the gradual emergence of numerous research papers and monographs on multiphase flow theory. It has also progressed from low concentration (dilute particle phase) flows to high concentration (dense particle phase) flows. Building upon the foundation of dilute phase flow, researchers have developed physical models and fundamental equations for dense phase flow\cite{11,12,13,14,15}. Subsequently, many researchers, including Campbell et al.\cite{16}, conducted a series of studies and established the particle dynamics. Then Ding and Gidaspow\cite{17} further developed and refined the theory of particle dynamics, introducing the concept of virtual temperature to characterize the effects of particle collisions in dense granular flows. They also conducted numerical simulations which yielded results consistent with experimental observations.

In the recent years, the complex interactions in multiphase flow between phases have garnered increasing attention. Many scholars have conducted detailed research on turbulent disturbances in two-phase flows, with particular emphasis on simulating gas-solid two-phase turbulence, which has become a research focus. Various approaches have been developed, including two-phase turbulence direct numerical simulation\cite{18}, Reynolds-averaged numerical simulation and large-eddy simulation\cite{lixianggui,eddy}.

Under extreme shock loading, particularly in two-phase problems, metal materials' surface probably experiences fragmentation, resulting in high-speed particles separating from the metal materials. This particulate matter then mixes with the surrounding gas, leading to a phenomenon called 'ejecat mixing'\cite{22,23}. The theoretical and experimental method and numerical simulation of ejecta mixing is currently a challenging and important issue in the field of implosion compression science \cite{18,19,20,21}. It represents one of the significant challenges in two-phase flow and is also a highly important area of application.

There are many experimental results publishing on the ejecta mixing. In 1983, Elias et al. published the experimental results of jet mixing involving metallic tin and argon gas \cite{24}. They conducted joint observations of ejecta mixing evolution using X-ray imaging and high-speed photography. Later Ogorodnikov et al. designed ejecta mixing experiments involving metallic lead and air\cite{25}. In the same year, Buttler et al. \cite{27} utilized PDV technology to investigate the velocity evolution of tin ejecta in neon gas. Their experimental results reflected the gradual decay of the ejecta over time due to fragmentation and resistance effects. In \cite{26}, Oro et al. conducted ejecta mixing experiments by preloading tungsten particles with a certain size distribution on an aluminum substrate using explosive loading, including vacuum ejecta mixing experiments, argon gas filling experiments, and hydrogen gas filling experiments. Sorenson et al. \cite{28} further studied the characteristics of ejecta size distribution with different gases.

On the other hand, studying particle motion through numerical simulation is also an important role to supplement experiments.In 2004, Wang et al.\cite{SPH} used SPH method to simulate the ejecta phenomenon on the metal surface under shock loading. Harrison, Fung, and others\cite{32}, based on the FLAG Lagrangian fluid dynamics program form, have developed an initial numerical simulation program for ejecta mixing. This program consists of two main models: the particle source model and the particle transport model. The particle source model is used to generate and initialize particles, including determining the conditions for ejecta initiation, namely the mass, the velocity distribution, the size distribution, the thermodynamic state, initial position of the particles and more. The particle transport model is given to simulate the motion of particles in the gas environment, with the main algorithm utilizing the MP-PIC method \cite{36,37,38}. In 2013, Fung et al. \cite{32} provided a detailed introduction to the settings of particles' initial state, but did not present quantitative calculations of ejecta mixing results or comparisons with experimental data.

Because of the complexity of the coupled gas-particle ejecta system, it means that while the fluid influences the particles, the particles also exert an effect on the fluid, necessiating the use of robust numerical methods for fluid motion simulations. In earlier years, a classical numerical method for fluid equation is primarily based on a conserved cell-centered grid discretization. The one-dimensional scheme of this method was initially presented by Godunov et al. \cite{39}, and later it was extended to multi-dimensional scenarios \cite{41,42}. Many scholars have proposed approximate Riemann solvers at the nodes, such as Despr{\'e}s and Mazera \cite{47}, Maire et al.\cite{maire2007cell}, which have shown significant improvements in correcting results.

Drawing from the advancements outlined in the above-mentioned literature, we employ a coupled gas-particle ejecta transporting model through configuring the initial states of the particles in a case of low concentration. And we focus on solving this problem using a cell-centered Arbitrary Lagrangian-Eulerian (ALE) method which is an extension of the method in \cite{shen2014robust2} from computational single-phase flow to the calculation of two-phase flow problems.

Based on the ejecta mixing and a coupled transporting effect, we apply the particle trajectory model to analyze the interactions and energy exchange between gas and particle phase. Through the numerical tests, we show the motion results of particles and fluid and a following effect of particles in fluid. Furthermore, we analyze the parameters that influence this effect and explore the reasons behind its formation.

The rest of the paper are organized as follows. In Section~2, we introduce the particle trajectory model and interactions between two phases. In Section~3, a numerical simulation method is given. In Section~4, we present the time discretization and a particle searching strategy for a CFL-like condition to limit the time step. In Section~5, we show the numerical tests, analyze the results and give a conclusion.

	
\section{Particle Trajectory Model}~\label{sec:main results}
\ 
\newline

In this section, we introduce the two-dimensional particle trajectory model by analyzing the dynamics of particles and compressible fluid. 

\subsection{Problem statement, Notation and Assumptions}\label{section:notation,assumptionm,problem}
\ 
\newline

Consider the following problem in a region $\Omega$ of two-dimensional space that multiple particles are mixed in the fluid(the total volume is less than a certain percentage). Regard each particle as an individual, interacting with the fluid separately and explore the motion of each particle independently. Since our particle concentration is less than a certain ratio, the particle flow is considered to be sparse particle flow, so it is assumed that each particle does not affect each other. And the fluid has initial density, velocity, pressure and interacts with particles to exchange momentum and energy at the same time.

In the following part, subscript $p$ indicates the physical parameter that is about particle phase, subscript $g$ is about the gas phase. Let $\rho_{g}$, $\boldsymbol{v_{g}}=(v_{g_x},v_{g_y})$, $E_{g}$, $\mu_g $, $T_g$, $C_g$ be the density, velocity, energy, viscosity,  thermal conductivity and $ P $ be the pressure of the fluid. $\boldsymbol{u_{p}}=(u_{p_x},u_{p_y})$, $\rho_{p}$, $\boldsymbol{v_{p}}=(v_{p_x},v_{p_y})$, $T_{p}$, $c_{p}$, $r_{p}$ denote the coordinate, density, velocity ,temperature, specific heat capacity and radius of particle. In particular, if there are multiple particles in the region $\Omega$ with different physical quantities, we can additionally use a subscript $m$ to identify them. Namely $\boldsymbol{u_{p_m}}$, $\rho_{p_{m}}$, $\boldsymbol{v_{p_{m}}}$, $T_{p_m}$, $c_{p_m}$ $r_{p_m}$denote the variable of the different particles, in which $m=1, \ldots, m(\Omega)$, $m(\Omega)$ is the total number of particles in the region $\Omega$.

To investigate the interactions between two phases, the force analysis of the particles in the system is essential. During the motion, particles are subject to gravity, drag, differential pressure, Basset force, Magnus force, Saffman force, etc. In general, not all of these forces are equally important. It is necessary to compare the magnitude of the above forces in order to estimate the relative importance of the various forces. Under normal conditions, other forces are several orders of magnitude smaller than the drag. So for the convenience of analysis, above forces except drag are not considered in the following model. 

In addition to the force analysis, the exchange of heat between the particles and the fluid is also significant. Normally there are three ways for transferring the heat including heat conduction, convection heat transfer, thermal radiation. Similarly, because of the orders of the magnitude, we only consider the convection heat transfer. 

To sum up, give the following assumptions:

\begin{itemize}
			\item[(1)]For each case we consider ideal compressible gas with ideal gas equation of state $P=(\gamma-1)\rho\epsilon$, where $\gamma$ is the adabatic index.
			\item[(2)]Viscosity and heat transfer exist only in interface interactions.
			\item[(3)]The motion of particles only take into account the drag due to friction.
            \item[(4)]The particles are ideal round or spherical, uniform in diameter.
            \item[(5)]There is no collision and interaction between particles.
            \item[(6)]All the energy produced by friction is absorbed by the fluid.
\end{itemize}

\subsection{Governing Equations}\label{section:equation}
\ 
\newline

According to the above assumptions, we introduce the equation of the system.

The equations for two-dimensional compressible particle trajectory model are as follows, for the fluid phase:


\begin{equation}\label{eq:fluid equation}
\begin{aligned}
\frac{\partial \mathbf{U}}{\partial t}+\frac{\partial \mathbf{F}(\mathbf{U})}{\partial x}+\frac{\partial \mathbf{G}(\mathbf{U})}{\partial y}=\mathbf{S}_{p},
\end{aligned}
\end{equation}\\
where the state vector and flux vectors are

\begin{equation}
	\mathbf{U}=\left[\begin{array}{c}
		\rho_g \\
		\rho_g v_{g_x} \\
		\rho_g v_{g_y} \\
		\rho_g E_g
	\end{array}\right], \quad \mathbf{F}(\mathbf{U})=\left[\begin{array}{c}
		\rho_g  v_{g_x} \\
		P + \rho_g v_{g_x}^2 \\
		\rho_g v_{g_x} v_{g_y} \\
		\rho_g E_g v_{g_x}+P v_{g_x}
	\end{array}\right], \quad \mathbf{G}(\mathbf{U})=\left[\begin{array}{c}
		\rho_g  v_{g_y} \\
		\rho_g v_{g_x} v_{g_y} \\
		\rho_g v_{g_y}^2+P \\
		\rho_g E_g v_{g_y}+P v_{g_y}
	\end{array}\right].
\end{equation}

For the particle phase, each particle has:
\begin{equation}\label{eq:particle equation}
\begin{aligned}
& \frac{\partial \boldsymbol{u_{p_m}}}{\partial t}=\boldsymbol{v}_{p_m}, \\
& \frac{\partial \boldsymbol{v}_{p_m} }{\partial t}=\boldsymbol{F}_{p_m}, \\ 
& \frac{4}{3} \pi \, r_p^3 \, \rho_p \, c_p \frac{\partial T_{p_m}}{\partial t}=Q_{p_m}, \\
\end{aligned}
\end{equation}\\
where the $\boldsymbol{F}_{p_m}$ and ${Q_{p_m}}$ indicate the interactions between the particle phase and fluid phase. 

And $\mathbf{S}_p$ is the resource term between the flow and particles

\begin{equation}
	\mathbf{S}_p=\left[\begin{array}{c}
		0 \\
		{F}_{p_x} \\
		{F}_{p_y} \\
		{E}_{p}
	\end{array}\right],
\end{equation}
where $\boldsymbol{F}_{p}$ is the force which is the source item of the particle's motion. ${E}_{p}$ is the exchange of energy including the heat and the kinetic energy.

In the particle trajectory model, the motion of every particles will be calculated respectively, actually we have

\begin{equation}
\boldsymbol{F}_{p}=\left[\begin{array}{c}
	{F}_{p_x} \\
	{F}_{p_y} \\
\end{array}\right]=\sum_{m=1}^{m(\Omega)}\boldsymbol{F}_{p_m},\\
\end{equation}
\begin{equation}
{E_{p}}=\sum_{m=1}^{m(\Omega)}{E}_{p_m}=\sum_{m=1}^{m(\Omega)}\boldsymbol{F}_{p_m} \cdot \boldsymbol{v}_{p_m}+Q_{p_m}.\\
\end{equation}

Then give the specific definition of interaction items $\boldsymbol{F}_{p_m}$ and ${Q_{p_m}}$. The force $\boldsymbol{F}_{p_m}$ comes from the drag due to friction. There are two drag models usually used during computing which are the Stoke's model and Crowe's model.

For the Stoke's model\cite{stokes2010richtmyer}, the following formulation is applied: 

\begin{equation}
\boldsymbol{F}_{p_m}=6 \pi r_p \mu_g (\boldsymbol{v}_g-\boldsymbol{v}_{p_m}),
\end{equation}
where $r_p$ is the radius of the particles and $\mu_g$ is the gas viscosity. 

For the Crowe's model\cite{crowe2011multiphase}, we have 

\begin{equation}
\boldsymbol{F}_{p_m}=\frac{1}{2} C_{d} K_{p_m} \rho_g \left|\boldsymbol{v}_g-\boldsymbol{v}_{p_m}\right|\left(\boldsymbol{v}_g-\boldsymbol{v}_{p_m}\right),
\end{equation}
where $C_{D}$ is the drag coefficient and $K_{p_m}$ is the upwind area. Drag coefficient $C_{D}$ is related to the Reynolds number $Re$, the following definition is applied\cite{morsi1972investigation}\cite{wang2012numericalCD}:

\begin{equation}\label{fdrag}
C_D=\frac{24}{Re} f_{d r a g}(Re).
\end{equation}

In formulation \eqref{fdrag}, $f_{d r a g}(Re)$ depends on the Reynolds number $Re$:

\begin{equation*}
R e=\frac{2 \rho_g r_{p_m}\left|\boldsymbol{v}_g-\boldsymbol{v}_{p_m}\right|}{\mu_g},
\end{equation*}
and

\begin{equation}
C_{\mathrm{D}}= 
\begin{cases}
\frac{24}{Re}, & Re<0.2, \\
\frac{24}{Re}\left(1+0.15 Re^{0.687}\right), & 0.2 \leqslant Re \leqslant 800, \\ 
0.5, & Re>800.\\
\end{cases}
\end{equation}

The calculation method of $C_{\mathrm{D}}$ is not same in different research. The above definition is given by Wang et al in \cite{wang2012numericalCD}.

$Q_{p_m}$ comes from the convection heat transfer\cite{chen2016heatq}:

\begin{equation}
Q_{p_m} = \frac{2 \mu_g}{P r} c_{pg} \pi r_{p_m}\left(T_{g}-T_{p_m}\right) Nu_{p_m},
\end{equation}
where $C_{pg}$ is constant-pressure specific heat capacity, $Nu_{p_m}$ is Nusselt numbers, $P r$ is Prandtl numbers of flow, of which specific definitions are as follows:

\begin{equation}
Pr=\frac{4 \gamma}{9 \gamma-5},
\end{equation} 
and

\begin{equation}
Nu_{p_m}=2+0.6 Pr^{1 / 3} Re^{1 / 2}.
\end{equation}

Based on the equation \eqref{eq:fluid equation} and \eqref{eq:particle equation}, we need to prove that the above gas-particle two phase system satisfies conservation law. It is obvious that the conservation of momentum because the drag between the particle and the fluid is a pair of interacting forces, namely the right items of the second equation of \eqref{eq:fluid equation} and of the equation \eqref{eq:particle equation}. The item $\boldsymbol{F}_{p_m} \boldsymbol{v}_{p_m}$ is the energy of the fluid which is converted into the kinetic energy of the particle and the convection heat transfer energy $Q_{p_m}$ is generated by the temperature difference.


\section{Numerical Scheme}\label{section:numerical scheme}
\ 
\newline

Then we give a numerical scheme to solve the equation \eqref{eq:fluid equation} and \eqref{eq:particle equation}.

\subsection{ALE scheme}\label{section:ALE}
\ 
\newline

For the subsequent discretization process, it is convenient to express equation \eqref{eq:fluid equation} in the form of a moving control volume formulation:

\begin{equation}\label{eq:control volumn formulation}
\frac{d}{d t} \int_{\Omega} \mathbf{U} d x d y+\int_{\partial \Omega}[(\mathbf{F}, \mathbf{G}) \mathbf{N}-(\mathbf{w} \cdot \mathbf{N}) \mathbf{U}] d l=\int_{\Omega} \mathbf{S}_{p} d x d y,
\end{equation}
$\mathbf{w}$ denotes the moving velocity of the control volume $\Omega$ while $\mathbf{N}$ denotes the unit outward normal direction along the boundary of $\Omega$. When $\mathbf{w}=\mathbf{u}$, the system simplifies into a Lagrangian formulation, and when $\mathbf{w}=\mathbf{0}$, it takes on an Eulerian form.

\begin{remark}
	It is important to note that in the equation \eqref{eq:control volumn formulation} the integration domain $\Omega$ is related to time $t$, that is $\Omega=\Omega \left( t \right)$. For convenience, we denote it as $\Omega$.
\end{remark}

\subsubsection{Notations on a generic polygonal grid}
\ 
\newline

\begin{figure}[htbp]
                \centering
                \includegraphics[width=0.3\textwidth]{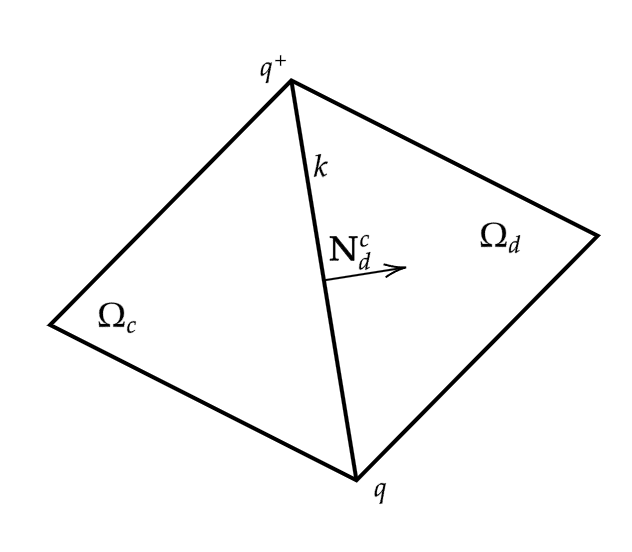}
            \caption{Notations on a grid}
            \label{fig:grid}    
\end{figure}

Every cell $\Omega_c$ within the mesh is assigned a distinct index represented by $c$. To refer to an adjacent cell $\Omega_d$ that shares a common edge with $\Omega_c$, we employ the subscript $d$. The edge where two units intersect is defined as $c \cap d$ or denoted as $k$, as shown in Fig \ref{fig:grid}. ${F}(c)$ denotes the set of adjacent cells surrounding the cell $\Omega_c$. ${Q}(c)$ is the set of all vertices of a cell $\Omega_c$. Each vertex within the mesh is allocated the index $q$, and $C(q)$ and $K(q)$ are defined as the sets of cells and edges respectively surrounding the vertex $q$, i.e.,

$$
\begin{aligned}
& C(q)=\left\{\Omega_c \text { : cells surrounding the vertex } q\right\}, \\
& K(q)=\{k \text { : edges surrounding the vertex } q\} .
\end{aligned}
$$

The physical quantities of fluid for example the density $\rho_{g_c}$, pressure $P_{g_c}$, velocity $\boldsymbol{v}_{g_c}$, energy $E_{g_c}, e_{g_c}$ are defined in the center of $\Omega_c$. The grid moving velocity, denoted as $\mathbf{w}_q$, is defined at the node $q$. Let $\mathbf{N}_d^c=\left(\left(n_x\right)_d^c,\left(n_y\right)_d^c\right)$ represents the unit normal vector of cell $\Omega_c$ along edge $c \cap d$, with the subscript $d$ and the superscript $c$ indicating that the vector direction originates from $\Omega_c$ and extends towards $\Omega_d$. 

Based on the edge's normal direction between two adjacent cells, we utilize indices $L$ and $R$ to represent the states on the left and right sides of an edge, respectively. For instance, $\mathbf{U}_L$ and $\mathbf{U}_R$ are state vectors in the left and right cells of the edge $c \cap d$.

To compute integration of the variable $\mathbf{F}$ across edge $\left[M_a, M_{a+1}\right]$, denote by $\mathbf{F}_{a,a+1}$ the fluxes on the face, $L_{a,a+1}$ is the length of $\left[M_a, M_{a+1}\right]$, we have

\begin{equation}
L_{a, a+1} \mathbf{F}_{a,a+1} \cdot \mathbf{N}_{\boldsymbol{a}, \boldsymbol{a}+\mathbf{1}}=\int_{M_a}^{M_{a+1}} \mathbf{F} \cdot \mathbf{N} \mathrm{d} l.
\end{equation}

Using the above transformation, we discretize the equations \eqref{eq:control volumn formulation} by utilizing the Godunov scheme on an unstructured mesh. Define cell-centered state variable $\mathbf{U}_c$ and resource term $\mathbf{S}_{p_c}$, then we have:

\begin{equation}\label{eq:dis control volumn}
\begin{aligned}
\frac{d}{d t} {\left|\Omega_c\right|} \mathbf{U}_c & ={\left|\Omega_c\right|} \mathbf{S}_{p_c}- \sum_d \int_d\left(\mathcal{T}_d^c\right)^{-1}\left[\mathbf{F}\left(\mathcal{T}_d^c \mathbf{U}\right)-\left(\mathbf{w}_d^c \cdot \mathbf{N}_d^c\right) \mathcal{T}_d^c \mathbf{U}\right] d l \\
& ={\left|\Omega_c\right|} \mathbf{S}_{p_c}- \sum_d L_d^c\left(\mathcal{T}_d^c\right)^{-1} \mathbf{F}_d^c\left(\mathbf{w}_d^c \cdot \mathbf{N}_d^c, \mathcal{T}_d^c \mathbf{U}_c, \mathcal{T}_d^c \mathbf{U}_d\right),
\end{aligned}
\end{equation}
where $\left|\Omega_c\right|$ represent the volume of cell $\Omega_c$. $\mathbf{w}_d^c$ is the moving velocity defined at the center of  cell edge $c \cap d$. $L_d^c$ is the length of the edge $c \cap d$. $\mathcal{T}_d^c$ and $\left(\mathcal{T}_d^c\right)^{-1}$ are the rotation matrix and its inverse as follows:

\begin{equation}
\mathcal{T}_d^c=\left(\begin{array}{cccc}
1 & 0 & 0 & 0 \\
0 & \left(n_x\right)_d^c & \left(n_y\right)_d^c & 0 \\
0 & -\left(n_y\right)_d^c & \left(n_x\right)_d^c & 0 \\
0 & 0 & 0 & 1
\end{array}\right), \quad\left(\mathcal{T}_d^c\right)^{-1}=\left(\begin{array}{cccc}
1 & 0 & 0 & 0 \\
0 & \left(n_x\right)_d^c & -\left(n_y\right)_d^c & 0 \\
0 & \left(n_y\right)_d^c & \left(n_x\right)_d^c & 0 \\
0 & 0 & 0 & 1
\end{array}\right) .
\end{equation}

There are many methods to calculate the approximate numerical flux  $\mathbf{F}_d^c$ between the cell $\Omega_c$ and $\Omega_d$. The most commonly used involves solving a one-dimensional Riemann problem along the normal direction of edge $c \cap d$ :
\begin{equation}\label{eq:Riemann problem}
\begin{aligned}
& \mathbf{U}_t+\mathbf{F}_x=\mathbf{S}_p, \\
& \mathbf{U}(x, 0)= \begin{cases}\mathbf{U}_L=\mathcal{T}_d^c \mathbf{U}_c, & x<0, \\
\mathbf{U}_R=\mathcal{T}_d^c \mathbf{U}_d, & x>0,\end{cases}
\end{aligned}
\end{equation}
and the numerical flux can be represented as follows:

\begin{equation}\label{eq: flux definition}
\mathbf{F}_d^c\left(w, \mathbf{U}_L, \mathbf{U}_R\right)=\mathbf{F}^w-w \mathbf{U}^w,
\end{equation}
where $w=\mathbf{w}_d^c \cdot \mathbf{N}_d^c$ is the normal projection of the moving velocity $\mathbf{w}_d^c$ , $\mathbf{U}^w$ and $\mathbf{F}^w$ are the projection of the state and flux functions onto the direction $x / t=w$ respectively.

\subsubsection{Flux formula of HLLC-2D on moving mesh}\label{sec:HLLC-2D}
\ 
\newline

We use the two-dimensional Riemann solver HLLC-2D method proposed in \cite{shen2014robust2} to evaluate the numerical flux. The construction of this method hinges on two important points. Firstly, it involves a nodal solver that relies on the vertex $q$ with state vectors $\left\{\mathbf{U}_c, c \in C(q)\right\}$ assigned to neighbor cells. Then the nodal contact velocity $\mathbf{u}_q^*$ is computed through this solver. Secondly, it involves the utilization of the HLLC-2D solver, an expanded scheme of the classical one-dimensional HLLC Riemann solver. In Riemann problem \eqref{eq:Riemann problem}, we have the nodal contact velocity $\mathbf{u}_q^*$, as well as the state vectors $\mathbf{U}_c$ and $\mathbf{U}_d$, which are on the sides of edge $c \cap d$. Through solving the Riemann problem, we could get the approximate numerical flux $\mathbf{F}_c^d$ and $\mathbf{F}_d^c$ and further solving the equation \eqref{eq:dis control volumn}.

In the following evaluation, assume that the nodal contact velocity $\mathbf{u}_q^*$ has already been calculated, in alignment with the preceding chapter's description. From the equation \eqref{eq:Riemann problem}, the initial values of the equation are:
\begin{equation*}
\mathbf{U}_L=\mathcal{T}_d^c \mathbf{U}_c, \quad \mathbf{U}_R=\mathcal{T}_d^c \mathbf{U}_d,
\end{equation*}
and the fluxes are $\mathbf{F}_L=\mathbf{F}\left(\mathbf{U}_L\right), \mathbf{F}_R=\mathbf{F}\left(\mathbf{U}_R\right)$.

In contrast to the classical one-dimensional Riemann solver, HLLC-2D method deviates by introducing an artificial contact wave speed denoted as $S_*=\mathbf{u}_q^* \cdot \mathbf{N}_d^c$, which is consistent with the nodal velocity $\mathbf{u}_q^*$.

Establish two velocities $S_L$ and $S_R$ and suppose $S_L \leqslant S_* \leqslant S_R$. As a result, the HLLC-2D solver comprises 4 constant state regions separated by $S_L, S_*$ and $S_R$. Within these regions, the vector of conservation quantities is defined as follows:

\begin{equation}
\mathbf{U}^w= \begin{cases}\mathbf{U}_L, & \text { if } w \leqslant S_L, \\ \mathbf{U}_L^*, & \text { if } S_L<w \leqslant S_*, \\ \mathbf{U}_R^*, & \text { if } S_*<w \leqslant S_R, \\ \mathbf{U}_R, & \text { if } w>S_R,\end{cases}
\end{equation}
where $\mathbf{U}_L^*, \mathbf{U}_R^*$ are the same as the classical HLLC solver. Within each state region, we introduce the two approximate fluxes which is from the different direction of the edge $c \cap d$:

\begin{equation}\label{eq:flux struction}
\mathbf{F}_H^w= \begin{cases}\mathbf{F}_{H, 1}, & \text { if } w \leqslant S_L, \\ \mathbf{F}_{H, 2}, & \text { if } S_L<w \leqslant S_*, \\ \mathbf{F}_{H, 3}, & \text { if } S_*<w \leqslant S_R, \\ \mathbf{F}_{H, 4}, & \text { if } w>S_R,\end{cases}
\end{equation}
where $H=L, R$, and the above flux scheme ensure strictly satisfies the Rankine-Hugoniot conditions. 

According to the equation \eqref{eq: flux definition} and flux struction \eqref{eq:flux struction}. Numerical fluxes encircling the vertex q on the moving mesh are as follows:

\begin{equation}
\mathbf{F}_{d, q}^c\left(w, \mathbf{U}_L, \mathbf{U}_R\right)=\mathbf{F}_L^w-w \mathbf{U}^w= \begin{cases}\mathbf{U}_L\left(u_L-w\right)+\mathbf{D}_L, & \text { if } w \leqslant S_L, \\ \mathbf{U}_L^*\left(S_*-w\right)+\mathbf{D}_L^*, & \text { if } S_L<w \leqslant S_*, \\ \mathbf{U}_R^*\left(S_*-w\right)+\mathbf{D}_L^*, & \text { if } S_*<w \leqslant S_R, \\ \mathbf{U}_R\left(u_R-w\right)+\mathbf{D}_R-\mathbf{D}_R^*+\mathbf{D}_L^*, & \text { if } w>S_R,\end{cases}
\end{equation}
and

\begin{equation}
\mathbf{F}_{c, q}^d\left(w, \mathbf{U}_L, \mathbf{U}_R\right)=\mathbf{F}_R^w-w \mathbf{U}^w= \begin{cases}\mathbf{U}_L\left(u_L-w\right)+\mathbf{D}_L-\mathbf{D}_L^*+\mathbf{D}_R^*, & \text { if } w<S_L, \\ \mathbf{U}_L^*\left(S_*-w\right)+\mathbf{D}_R^*, & \text { if } S_L \leqslant w<S_*, \\ \mathbf{U}_R^*\left(S_*-w\right)+\mathbf{D}_R^*, & \text { if } S_* \leqslant w<S_R, \\ \mathbf{U}_R\left(u_R-w\right)+\mathbf{D}_R, & \text { if } w \geqslant S_R,\end{cases}
\end{equation}
where $\mathbf{D}_H=(0, P, 0, P u)^T, \mathbf{D}_H^*=\left(0, P_H^*, 0, P_H^* S_*\right)^T$ for $H=L, R$ are vectors including the pressure and work terms on ends and in the middle two regions, respectively. Pressures are

\begin{equation}\label{eq:pressure}
\left\{\begin{array}{l}
P_L^*=p_L+\rho_L\left(S_L-u_L\right)\left(S_*-u_L\right), \\
P_R^*=p_R+\rho_R\left(S_R-u_R\right)\left(S_*-u_R\right) .
\end{array}\right.
\end{equation}

An analogous derivation can be given at vertex $q^{+}$, and flux $\mathbf{F}_{d, q^{+}}^c$ is evaluated. For the edge $c \cap d$, the flux on it finally is:

\begin{equation}
\mathbf{F}_d^c\left(w, \mathcal{T}_d^c \mathbf{U}_c, \mathcal{T}_d^c \mathbf{U}_d\right)=\frac{1}{2}\left(\mathbf{F}_{d, q}^c\left(\mathbf{w}_q \cdot \mathbf{N}_d^c, \mathcal{T}_d^c \mathbf{U}_c, \mathcal{T}_d^c \mathbf{U}_d\right)+\mathbf{F}_{d, q^{+}}^c\left(\mathbf{w}_{q^{+}} \cdot \mathbf{N}_d^c, \mathcal{T}_d^c \mathbf{U}_c, \mathcal{T}_d^c \mathbf{U}_d\right)\right).
\end{equation}

\subsubsection{Nodal Solver}\label{subsec:node solver}
\ 
\newline

In Subsection \ref{sec:HLLC-2D}, we have assumed that the nodal contact velocity is already calculated. And since $\mathbf{F}_d^c \neq \mathbf{F}_c^d$ may happen, the conservation in cells is not satisfied. Therefore we introduce a nodal solver that consistently evaluates the $\mathbf{u}_q^*$ , ensuring the scheme preserves local conservation of mass, momentum and energy. This procedure can be derived from the global conservation within domain $\Omega$. By maintaining the local nodal conservation including mass, momentum and total energy, contact velocity could be calculated. By using the pressure equation \eqref{eq:pressure}, we have:

\begin{equation}\label{eq:nodal solver}
\sum_{k \in K(q)} L_k\left(\alpha_{L, k}+\alpha_{R, k}\right)\left[\mathbf{u}_q^* \cdot \mathbf{N}_k-v_k^*\right] \mathbf{N}_k=0,
\end{equation}
where

\begin{equation}
\alpha_{L, k}=-\rho_{L, k}\left(S_{L, k}-u_{L, k}\right), \quad \alpha_{R, k}=\rho_{R, k}\left(S_{R, k}-u_{R, k}\right),
\end{equation}
and $v_k^*$ represents the contact velocity in the classical one-dimensional HLLC Riemann solver for edge $k$

\begin{equation}
v_k^*=\frac{P_{L, k}-P_{R, k}+\alpha_{L, k} u_{L, k}+\alpha_{R, k} u_{R, k}}{\alpha_{L, k}+\alpha_{R, k}} .
\end{equation}

Solve the equation \eqref{eq:nodal solver} to calculate the contact velocity $\mathbf{u}_q^*$ at each vertex and we have

\begin{equation}\label{eq:nodal velocity}
\mathbf{u}_q^*=M^{-1} \sum_{k \in K(q)} L_k\left(\alpha_{L, k}+\alpha_{R, k}\right) v_k^* \mathbf{N}_k,
\end{equation}
where

$$
M=\left[\begin{array}{cc}
\sum_{k \in K(q)} L_k\left(\alpha_{L, k}+\alpha_{R, k}\right) n_{x, k}^2 & \sum_{k \in K(q)} L_k\left(\alpha_{L, k}+\alpha_{R, k}\right) n_{x, k} n_{y, k} \\
\sum_{k \in K(q)} L_k\left(\alpha_{L, k}+\alpha_{R, k}\right) n_{x, k} n_{y, k} & \sum_{k \in K(q)} L_k\left(\alpha_{L, k}+\alpha_{R, k}\right) n_{y, k}^2
\end{array}\right] .
$$

\subsubsection{Boundary Condition}
\ 
\newline

In this ALE method, the conventional zero-order extrapolation method is employed to enforce boundary conditions. For velocity-based boundaries, there are two different cases:

\begin{itemize}
	\item[(1)] When the three boundary nodes are not collinear, a straightforward extrapolation method is implemented.
	\item[(2)] When the three boundary nodes are collinear, we employ the mirror extrapolation method, which involves constructing a symmetrical point for each vertex connected to a boundary node.
\end{itemize}
The above method ensures the same scheme \eqref{eq:nodal velocity} as used for internal vertex when calculating $\mathbf{u}_q^*$.


\section{Time discretization and particle searching}\label{section:time discretization}
\ 
\newline
	
We discretize in time for the system \eqref{eq:dis control volumn} that describes the evolution of the variables in the cell $\Omega_i$. Suppose that there are $m(\Omega_i)$ particles in cell $\Omega_i$, of which the index is $m=1, \ldots, m(\Omega_i)$, then we use the second-order Runge-Kutta method to obtain the following fully discretized system:

\begin{equation}\label{eq:dis fluid1}
\begin{aligned}
\hat{\mathbf{U}}_{i}^{n+1} & =\frac{\Delta t^n}{\left|\Omega_i^{n+1}\right|} \sum_m \mathbf{S}_{p_i}^{n}+\frac{\left|\Omega_i^n\right|}{\left|\Omega_i^{n+1}\right|} \mathbf{U}_{i}^n-\frac{\Delta t^n}{\left|\Omega_i^{n+1}\right|} [\sum_d L_d^i\left(\mathcal{T}_d^i\right)^{-1} \mathbf{F}_d^i\left(\mathbf{w}_d^i \cdot \mathbf{N}_d^i, \mathcal{T}_d^i \mathbf{U}_i^{n}, \mathcal{T}_d^i \mathbf{U}_d^{n}\right)], \\
\end{aligned}
\end{equation}

\begin{equation}
	\overline{\mathbf{U}}_{i}^{n}=\frac{1}{2}({\mathbf{U}}_{i}^{n}+\hat{\mathbf{U}}_{i}^{n+1}),
\end{equation}

\begin{equation}\label{eq:dis fluid}
\begin{aligned}
\mathbf{U}_{i}^{n+1}=\frac{\Delta t^n}{\left|\Omega_i^{n+1}\right|} \sum_m\mathbf{S}_{p_i}^{n}+\frac{\left|\Omega_i^n\right|}{\left|\Omega_i^{n+1}\right|} \mathbf{U}_{i}^n-\frac{\Delta t^n}{\left|\Omega_i^{n+1}\right|} \sum_d L_d^i\left(\mathcal{T}_d^i\right)^{-1}
\mathbf{F}_d^i\left(\mathbf{w}_d^i \cdot \mathbf{N}_d^i, \mathcal{T}_d^i \overline{\mathbf{U}}_i^{n}, \mathcal{T}_d^i \overline{\mathbf{U}}_d^{n}\right),\\
\end{aligned}
\end{equation}
and similarly

\begin{equation}\label{eq:dis particle}
\begin{aligned}
& \hat{\boldsymbol{v}}_{p_m}^{n + 1}=\boldsymbol{v}_{p_m}^{n} + \frac{\Delta t^n}{\pi r_{p_m}^2 {\rho}_{p_m}} \, \boldsymbol{F}_{p_m}^{n}, \\
& \hat{\boldsymbol{F}}_{p_m}^{n+1}={\boldsymbol{F}}_{p_m}(\frac{1}{2}(\hat{\boldsymbol{v}}_{p_m}^{n + 1}+{\boldsymbol{v}}_{p_m}^{n})), \\ 
& \hat{T}_{p_m}^{n + 1}=T_{p_m}^{n}+{\Delta t^n}\, \frac{3}{4 \pi  r_{p_m}^3  \rho_{p_m}  c_{p_m}} Q_{p_m}^{n}, \\
& \hat{Q}_{p_m}^{n+1}={Q}_{p_m}(\frac{1}{2}({T}_{p_m}^{n}+\hat{T}_{p_m}^{n + 1})).
\end{aligned}
\end{equation}

so that 

\begin{equation}\label{eq:dis particle}
\begin{aligned}
& \boldsymbol{v}_{p_m}^{n + 1}=\boldsymbol{v}_{p_m}^{n} + \frac{\Delta t^n}{\pi r_{p_m}^2 {\rho}_{p_m}} \, \hat{\boldsymbol{F}}_{p_m}^{n+1}, \\
& \boldsymbol{u}_{p_m}^{n + 1}=\boldsymbol{u}_{p_m}^{n} + {\Delta t^n} \, \hat{\boldsymbol{v}}_{p_m}^{n + 1}, \\
& T_{p_m}^{n + 1}=T_{p_m}^{n}+\Delta t^n \, \frac{3}{4 \pi  r_{p_m}^3  \rho_{p_m}  c_{p_m}} \hat{Q}_{p_m}^{n+1}. \\
\end{aligned}
\end{equation}

\subsection{Particle searching}\label{section:particle searching}
\ 
\newline

Before computing the time step $\Delta t^n=t_{n+1}-t_n$, we need to give the process of searching the particles during a time step. In the numerical simulation, the interaction between the particles and flow is calculated in every cells instead of the whole space. That is, if consider the cell $\Omega_i$, the momentum and energy exchanged between the two phase is that all the particles in cell $\Omega_i$ interact with the fluid in $\Omega_i$. So before compute the physical properties we should know which particles are in the cell. At the beginning of time, the coordinates of particles are given. And after every time steps, we should search the particles and confirm the location of every particles. In the numerical simulations, we use the following method for searching the particles:

\begin{itemize}
			\item[(1)] Suppose that we know the coordinates of particles in the cell $\Omega_i$ at time $t_n$, namely $\boldsymbol{u_{m}^n}$, $m=1, \ldots, m(\Omega_i)$.
			\item[(2)] Calculate the equation \eqref{eq:dis particle} and get the coordinate of particle m, $\boldsymbol{u_{m}^{n+1}}$ , at time $t_{n+1}$.
			\item[(3)] Search the particle m in the updated geometrical characteristics of the cell. First determine whether the particle m is in the cell $\Omega_i$. If then, proceed to step 5.
            \item[(4)] If not, search the particle m in the cell $\Omega_c$ which is the surrounding cells of $\Omega_i$, $c=1,\dots,n(i)$.
            \item[(5)] Update the position of particle m. 
\end{itemize}

\begin{figure}[htbp]
                \centering
                \includegraphics[width=0.8\textwidth]{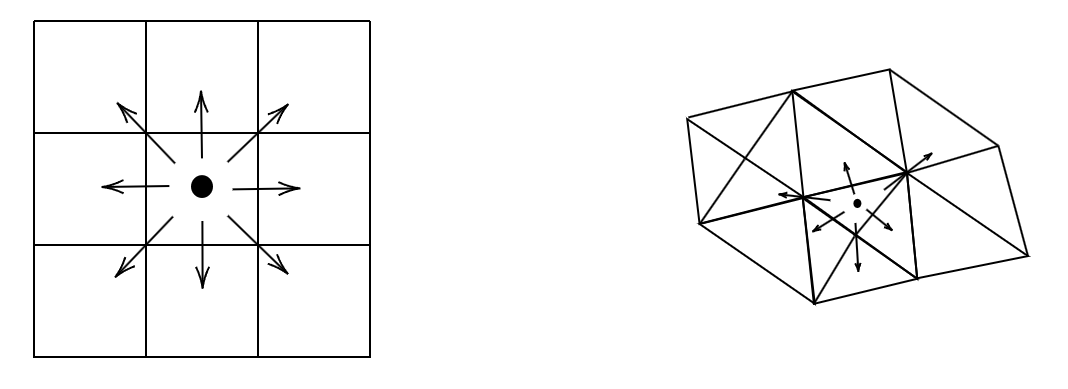}
            \caption{Particle searching}
            \label{fig:particle_searching}    
\end{figure}

Using the above method, we could get the updated position of particle k under the condition that we give a limitation to the time step $\Delta t^n$. In figure \ref{fig:particle_searching} we give the schematic diagrams of particle searching on structured and unstructured mesh respectively.

\subsection{Time step limitation}\label{section:time step limitation}
\ 
\newline

Because of using the explicit method, we need to observe some limitations on time step. They are as follows:

\begin{itemize}
			\item[(1)] A CFL condition to be a priori satisfied for entropy inequality.
			\item[(2)] The variation of cell volume corresponds coherently with the motion of the vertices.
			\item[(3)] Make sure that the particles can be searched.
\end{itemize}

Assume that the physical properties and geometric characteristics of the cell $\Omega_i$ are known at the time $t_n$. We will compute the variables' value at time $t_{n+1}$.

\subsubsection{CFL condition}\label{section:cfl condition}
\ 
\newline

First we derive a CFL condition to ensure  positive entropy production within cell $\Omega_i$ during the computation. At time $t_n$, we denote by $\lambda_i^n$ the inner circle radius of the cell $\Omega_i$. Then

\begin{equation}
\Delta t_e=C_e \min _{i=1, \ldots, I} \frac{\lambda_i^n}{\sqrt{(u_{g_i}^n)^2+(v_{g_i}^n)^2} + c_i^n},
\end{equation}
where $C_E$ is a strictly positive coefficient and $c_i$ is the sound speed of the cell $\Omega_i$. According to the numerical experiments, the value of $C_E$ between $0.3~0.4$ usually could get a good numerical results

\subsubsection{Criterion on the variation of volume}
\ 
\newline

We also establish a criterion for the variation of cell volume. At time $t_n$, we can get the contact velocity of vertex, we have:

\begin{equation}
\Delta t_v=C_V \min _{i=1, \ldots, I}\left\{\frac{\left| A_i^n \right|}{\left| A_i^{\prime} \right|}\right\} ,\\
\end{equation}
where $C_v$ can be set $0.1$ and 

\begin{equation}
\left|A_i^{\prime}\right|=\frac{1}{2} \sum_f L_f^c\left(\mathbf{u}_q^*+\mathbf{u}_{q^{+}}^*\right) \cdot \mathbf{N}_f^c .
\end{equation}

\subsubsection{Particle searching limitation}
\ 
\newline

Refer to Subsection \ref{section:particle searching}, suppose that particle $m$ is located in $\Omega_{i}$, we check the cell $\Omega_{i}$ and its surrounding cells $\Omega_{c}$, $c=1,\dots,n(i)$. From time $t_n$ to $t_{n+1}$ the particles should be specify the range of motion of the particle.  In the subsection \ref{section:cfl condition} denote by $\lambda_{k}^n$ the inner circle radius of the cell $\Omega_{i}$ and the cell $\Omega_{c}$. We have:

\begin{equation}
\Delta t_p=C_p \min_{m} \min _{k=1, \ldots, n(i)+1} \frac{\lambda_{k}^n}{|\boldsymbol{v_{p_m}^n}|},\\
\end{equation}
where $n(m)$ is the number of cells containing and surrounding the particle. In numerical simulation we take $C_p=0.5$.

Finally, the estimation of the next time step the estimation of the next time step $\Delta t^{n}=t_{n+1}-t_{n}$ is given by
\begin{equation}\label{eq:time step}
\Delta t^{n}=\min \left(\Delta t_e, \Delta t_v, \Delta t_p , 1.01*\Delta t^{n-1}\right).\\
\end{equation}

\subsection{Algorithm for full discretization system}
\ 
\newline

Give the description of the whole computation process:\\
\textbf{Step 1}. Initialization
\ 
\newline

At time $t=t^n$, in each cell $\Omega_i, \ i=1, \ldots, I$: there are particle variables: $\boldsymbol{u}_{p_m}^n$, $\rho_{p_m}^n$, $\boldsymbol{v}_{p_m}^n$, $T_{p_m}^n$, $r_{p_m}^n$, $m=1, \ldots, m(\Omega_i)$. Fluid variables: $\boldsymbol{v_{g}^n}$, $E_{g}^n$, $\rho_g^n $, $T_g^n$ and the geometrical quantities: $L_{k}^n$, $\mathbf{N}_{k}^n$ for $k \in K(q) $.\\
\ 
\newline
\textbf{Step 2}. Nodal solver of flow
\ 
\newline

The contact velocity $\boldsymbol{u}_q^*$ at vertex $q$ is calculated by solving the equations in subsection \ref{subsec:node solver}.\\
\ 
\newline
\textbf{Step 3}. Time step limitation
\ 
\newline

Obtain the time step $\Delta t^{n}$ according to the \eqref{eq:time step}.\\
\ 
\newline
\textbf{Step 4}. Grid moving strategy
\ 
\newline

In the ALE method, we need to give a specific grid-moving strategy, such as harmonic mapping equation, to compute vertex coordinate. And the nodal velocity is calculated by $\mathbf{w}_q=\left(\mathbf{x}_q^{n+1}-\mathbf{x}_q^n\right) / \Delta t^{n}$. \\
\ 
\newline
\textbf{Step 5}. Interaction between two phases
\ 
\newline

Calculate the interactions $\mathbf{F}_{p}$ and ${E_{p}}$ at time $t_n$.\\
\ 
\newline
\textbf{Step 6}. The Riemann solver HLLC-2D
\ 
\newline

At the edge $i \cap d$ of cell $\Omega_i$ and $\Omega_d$, compute the fluxes $\mathbf{F}_{d, q}^i$ and $\mathbf{F}_{d, q^{+}}^i$. The edge flux $\mathbf{F}_d^i$ is the average of $\mathbf{F}_{d, q}^i$ and $\mathbf{F}_{d, q^{+}}^i$.
\ 
\newline
\textbf{Step 7}. Update the state variable of flow
\ 
\newline

Update the physical quantities in cell $\Omega_i$. \\
\ 
\newline
\textbf{Step 8}. Update the particle
\ 
\newline

Calculate the variables of particles by equation \eqref{eq:dis particle} and search the position of particle at time $t_{n+1}$. Then returning to Step 1 if the calculation continues.


\section{Nemurical tests and results}\label{section:simulation and results}
\ 
\newline

In the following section, we present several examples to validate the above numerical scheme. We begin by several classical one-dimensional particle transport tests in quiescent fluid, in fluid under constant acceleration and in fluid under sinusoidal acceleration. The relationship between particle size and particle motion path is further explored in this part. Next we explore the motion of fluid in gas-particle Sod test. Finally, we give a two-dimensional practical test case about plumbum plate ejecta motion under the high speed shock. For generality, we calculate these numerical tests on unstructured mesh.

\subsection{Particle transport tests}
\ 
\newline

In this test, we investigate the motion of particles immersed in a fluid under different state which is an important problem in many applications.

\subsubsection{Transport in a quiescent fluid}
\ 
\newline

\begin{figure}[htbp]
              \begin{minipage}[t]{0.5\textwidth}
                \centering
                \includegraphics[width=\textwidth]{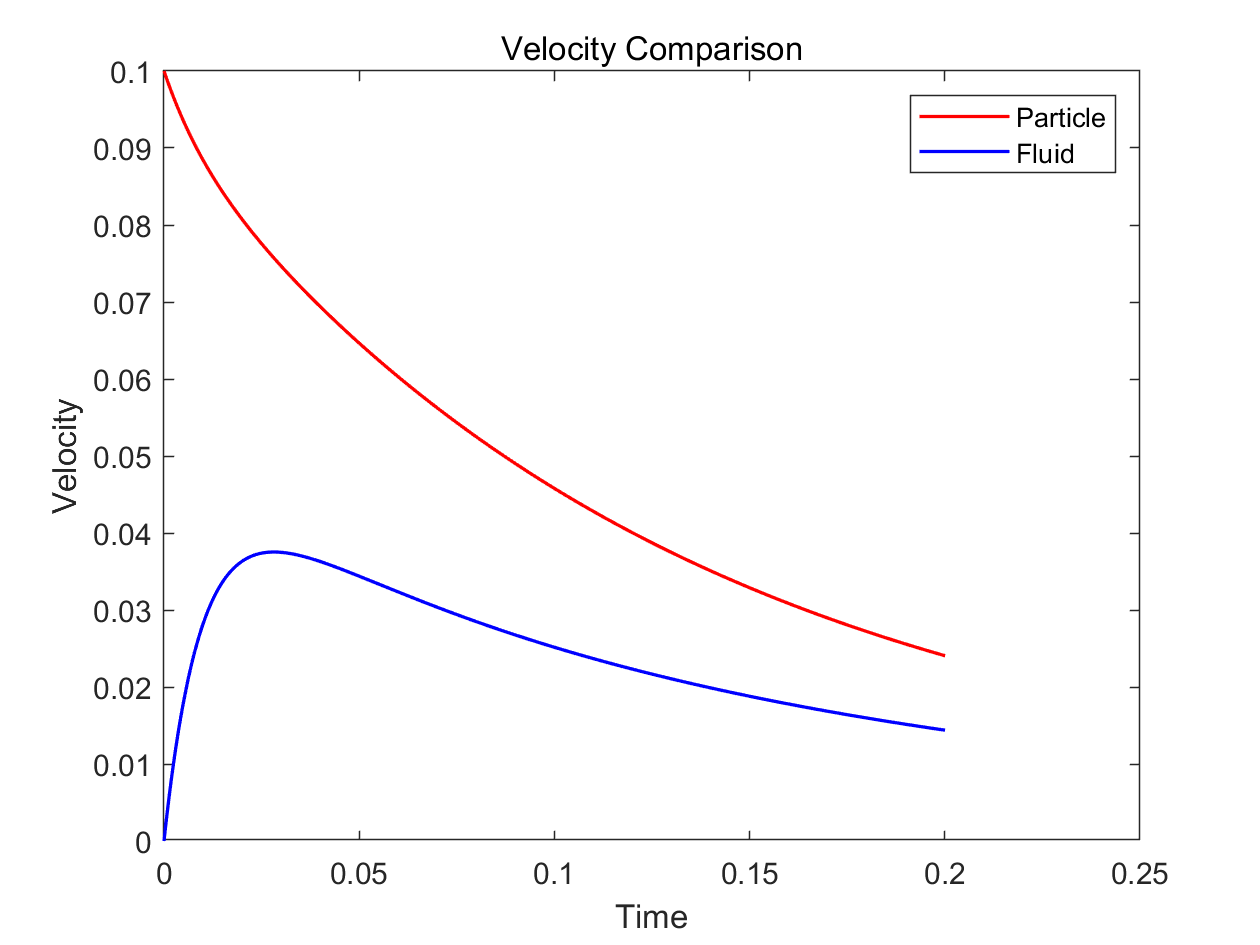}
                $t=0.2$
                \end{minipage}\hfill
                \begin{minipage}[t]{0.5\textwidth}
                 \centering
                \includegraphics[width=\textwidth]{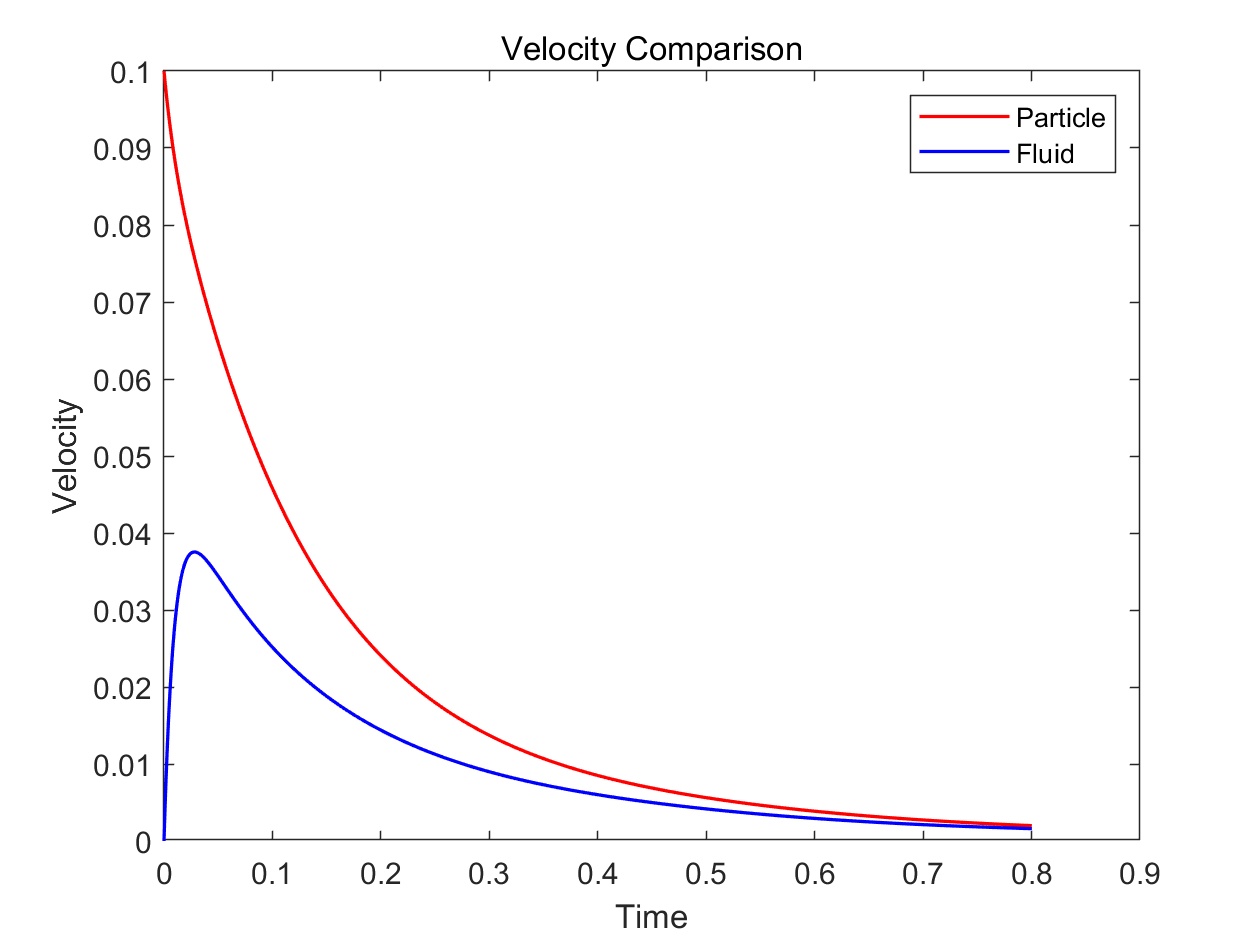}
                $t=0.8$
            \end{minipage}
            \caption{Drag model}   
            \label{fig:1}    
\end{figure}

The motion of a particle acted upon by no force except drag is the most basic model. We compute this problem to investigate the particle trajectories under the drag influence. Suppose a nonzero initial speed for the particle at the middle of a nozzle filled with fluid and compare the velocity of particle and fluid. According to the Figure \ref{fig:1}, the velocity of the particle is decreasing during the motion, with the deceleration gradually decreasing. At the same time, the speed of the fluid reaches the peak rapidly within a short period and then decreases. The speed difference between the particle and fluid ultimately converges, contrary to the intuition that the fluid velocity keeps increasing until they become consistent.

\begin{figure}[htbp]
              \begin{minipage}[t]{0.5\textwidth}
                \centering
                \includegraphics[width=\textwidth]{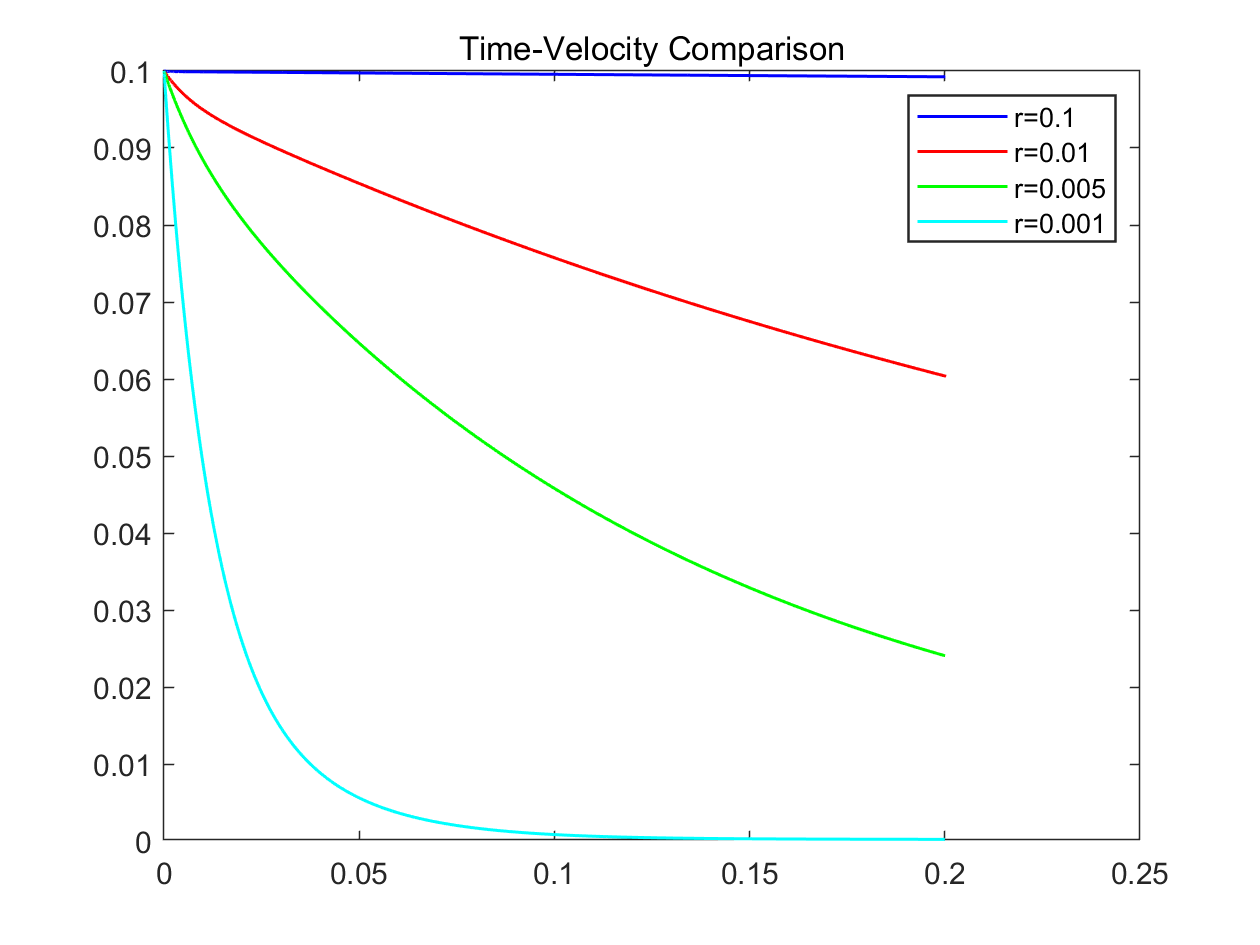}
                Time-Velocity
                \end{minipage}\hfill
                \begin{minipage}[t]{0.5\textwidth}
                 \centering
                \includegraphics[width=\textwidth]{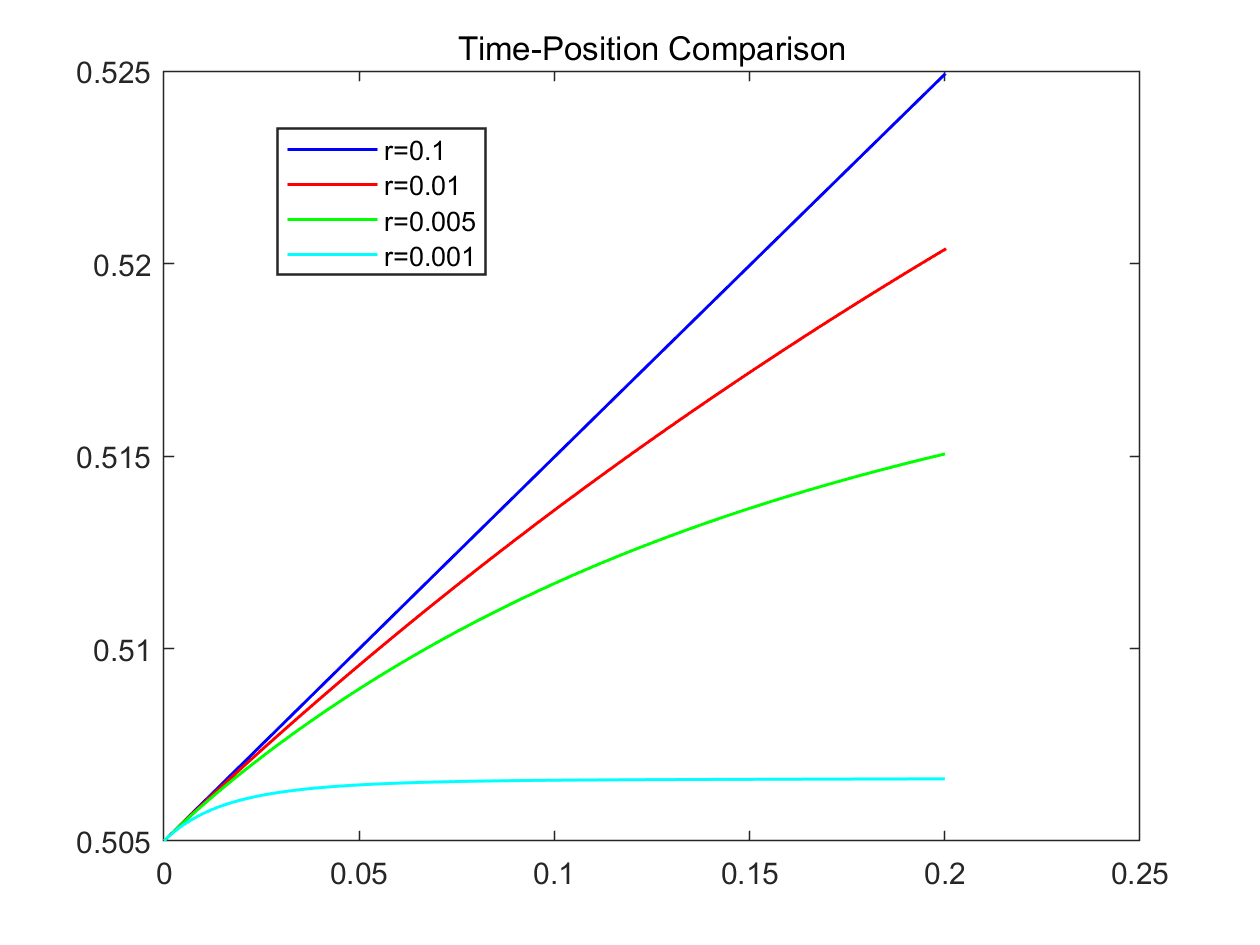}
                Time-Position
            \end{minipage}
            \caption{Different radius comparison in a quiescent fluid}   
            \label{fig:2}    
\end{figure}

Next based on the above tests, we further conduct the extend tests, investigating the speed variations by setting different radius for the particles in Figure \ref{fig:2}. Obviously, as the particle radius increases, the speed decreases more slowly. A noteworthy phenomenon is that when the particle becomes sufficiently large, it is almost unaffected by the fluid, a characteristic reflected in subsequent numerical results.

\subsubsection{Motion under constant acceleration}
\ 
\newline

\begin{figure}[htbp]
                \centering
                \includegraphics[width=1\textwidth]{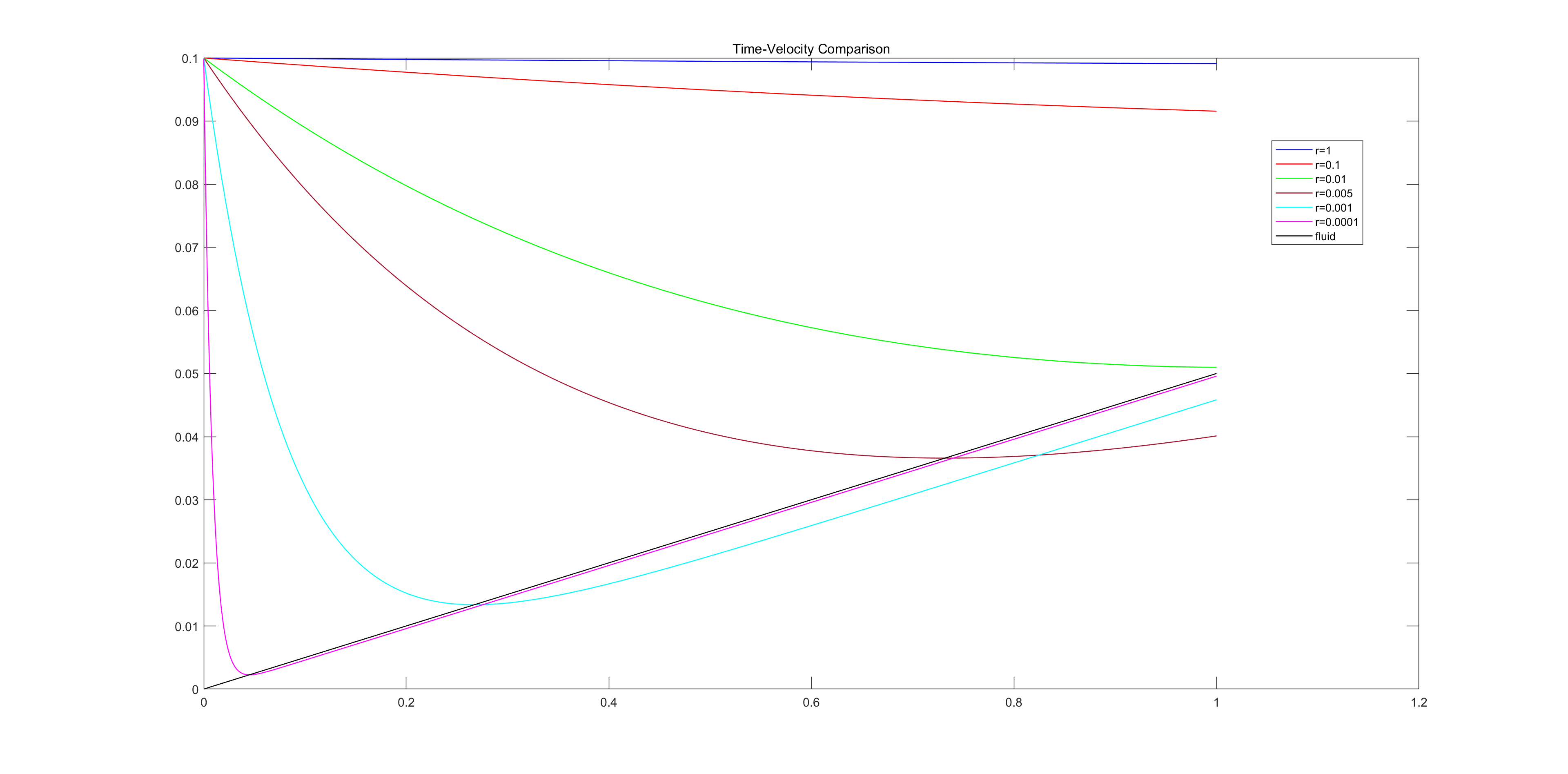}
                Time-Velocity
            \caption{Different radius comparison in a constant acceleration fluid}
            \label{fig:3.1}    
\end{figure}

\begin{figure}[htbp]
                \centering
                \includegraphics[width=1\textwidth]{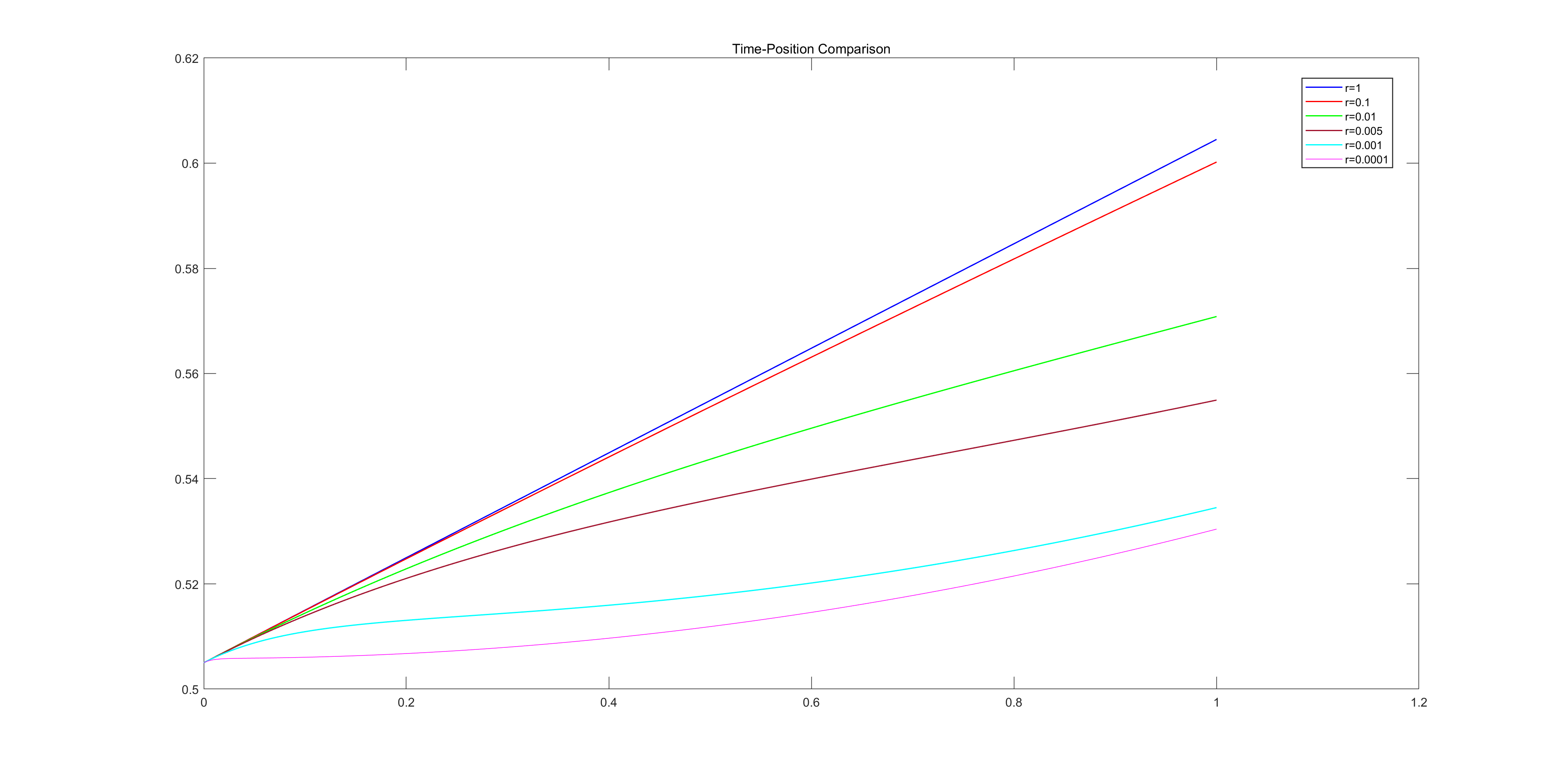}
                Time-Position
            \caption{Different radius comparison in a constant acceleration fluid}   
            \label{fig:3.2}    
\end{figure}

In the second series of tests, according to the Figure \ref{fig:3.1} and \ref{fig:3.2}, the different-size particles are given a nonzero initial velocity into a quiescent fluid with a constant acceleration. At the start of the motion, the fluid velocity is less than all the particles when the particles are enforced to decelerate, and then accelerate when the fluid velocity exceeds the particle velocity. Analyzing the curves of $r=0.0001$ and $r=0.001$, it can be observed that once the particle velocity is lower than that of the fluid, the particle accelerates with the same acceleration as the fluid, consistently lagging behind it. Similarly, when the particle becomes sufficiently large, for instance $r=1$, it is scarcely affected.

In addition, it is obvious that the particles eventually move with the different velocity differences from the fluid. Suppose that the acceleration of fluid is $\boldsymbol{a}_g$, the velocity difference is $\Delta\boldsymbol{v}=\boldsymbol{v}_g-\boldsymbol{v}_p$. There is the equation that $\boldsymbol{a}_g=\boldsymbol{F_d}/m_p=\boldsymbol{a}_p$ when the motion reaches the balance. In a general physical computation, drag force is positive correlation with the $\Delta\boldsymbol{v}$, we have $\boldsymbol{F_d}=\kappa \Delta\boldsymbol{v}$ in which $\kappa$ is a variable positive coefficient related to radius, density and viscosity when we used the different drag force model affecting the final velocity difference. For instance, in Stoke's drag model, there are $\kappa=6 \pi r_p \mu_g$ and $\Delta\boldsymbol{v}=\boldsymbol{a}_g r_p /6 \mu_g$ which is consistent with the numerical results.

\subsubsection{Motion under sinusoidal acceleration}
\ 
\newline

\begin{figure}[htbp]
              \begin{minipage}[t]{0.45\textwidth}
                \centering
                \includegraphics[width=\textwidth]{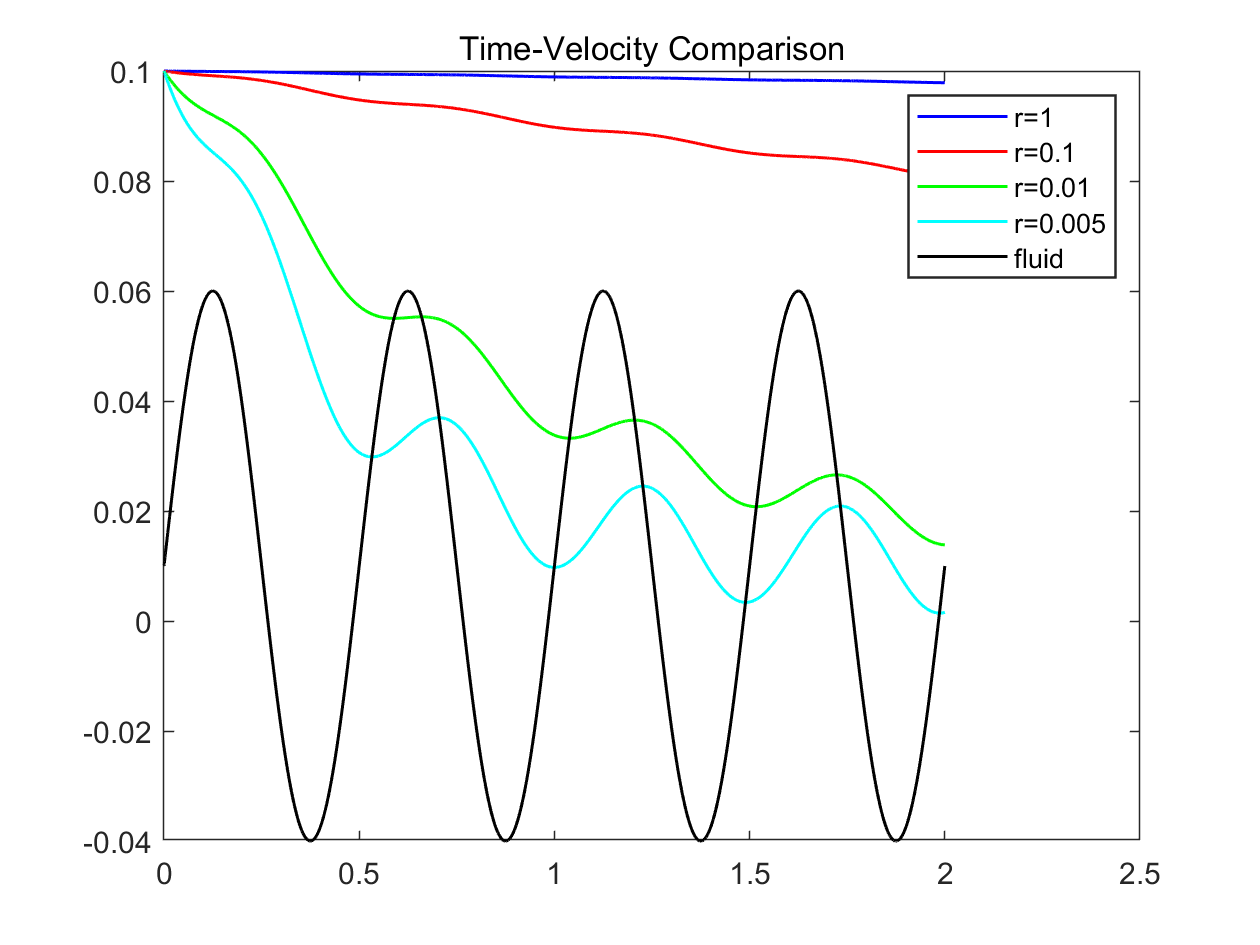}
                Time-Velocity
                \end{minipage}\hfill
                \begin{minipage}[t]{0.45\textwidth}
                 \centering
                \includegraphics[width=\textwidth]{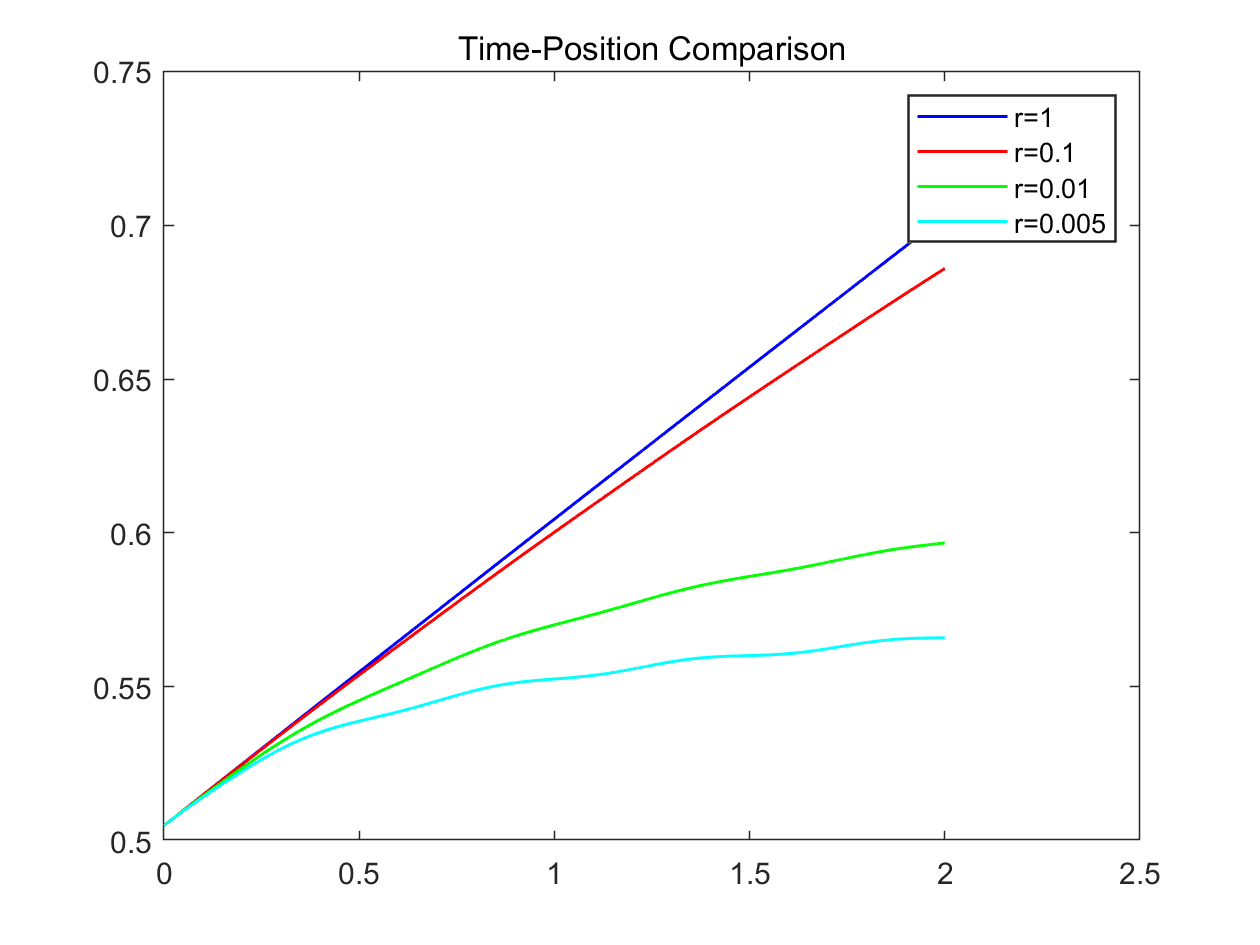}
                Time-Position
            \end{minipage}
            \caption{Different radius comparison in a sinusoidal acceleration fluid}   
            \label{fig:4}    
\end{figure}

In the final test,  the different-size particles are given a nonzero initial velocity into a quiescent fluid with a sinusoidal acceleration, with the result in Figure \ref{fig:4}. Comparing with the black line 

\subsection{Multiphase Sod test}
\ 
\newline

We are concerned with a shocktube problem\cite{sodtest} in the region $[-0.5,0.5] \times [0,0.1] $ of which solution consists of a left rarefaction, a contact and a right shock. We consider a shock tube of unity length. We test a series of initial conditions containing the particles of different size and analysis the influence of interaction to the flow. In the simulation, particles are uniformly distributed in the all cells of the region and the flow phase is an ideal gas. At beginning, the states on the left and right sides of $x=0$ are constant. The left state is a high pressure fluid set as $\left(\rho_l, P_l, u_l\right)=(9.6,1013250,0)$, and the right state is a low pressure fluid defined by $\left(\rho_r, P_r, u_r\right)=(1.2,101325,0)$, see \cite{tbl}. 

For the multiphase Sod shock tube, there is no exact solution. But we can analysis the physical interaction and results according to the solution of numerical simulation. In the numerical tests, we adjust the particles' radius which achieve the different volume percentages, namely $0\%$, $0.01\%$, $4\%$, to research the effects on the fluid motion. Then in the table \ref{tab:ic particles}, we have the initial conditions for particles.

\begin{table}\label{tab:ic particles}
\caption{\textbf{Initial conditions for particles}}
\centering
\begin{tabular}{cccc}
	\hline
		Case	   \quad&volume percentage \quad&radius  \quad& velocity \\
	\hline
		1	             \quad&  $0\%$   \quad&  $0$      \quad&  0\\

		2	             \quad&  $0.01\%$    \quad&  $0.0001$    \quad&  0  \\
		
		3				\quad&  $4\%$      \quad&  $0.002$        \quad&  0 \\
		
	\hline
	\end{tabular}	
\end{table}

\begin{figure}[htbp]
	\begin{minipage}[t]{0.5\textwidth}
		\centering
		\includegraphics[width=\textwidth]{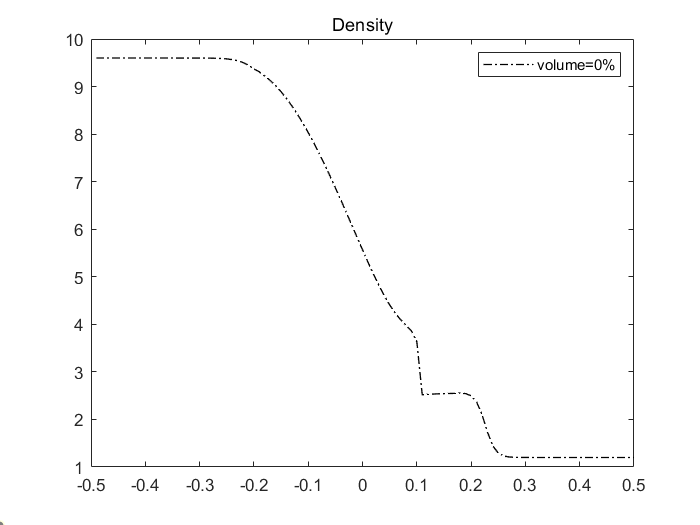}
	\end{minipage}\hfill
	\begin{minipage}[t]{0.45\textwidth}
		\centering
		\includegraphics[width=\textwidth]{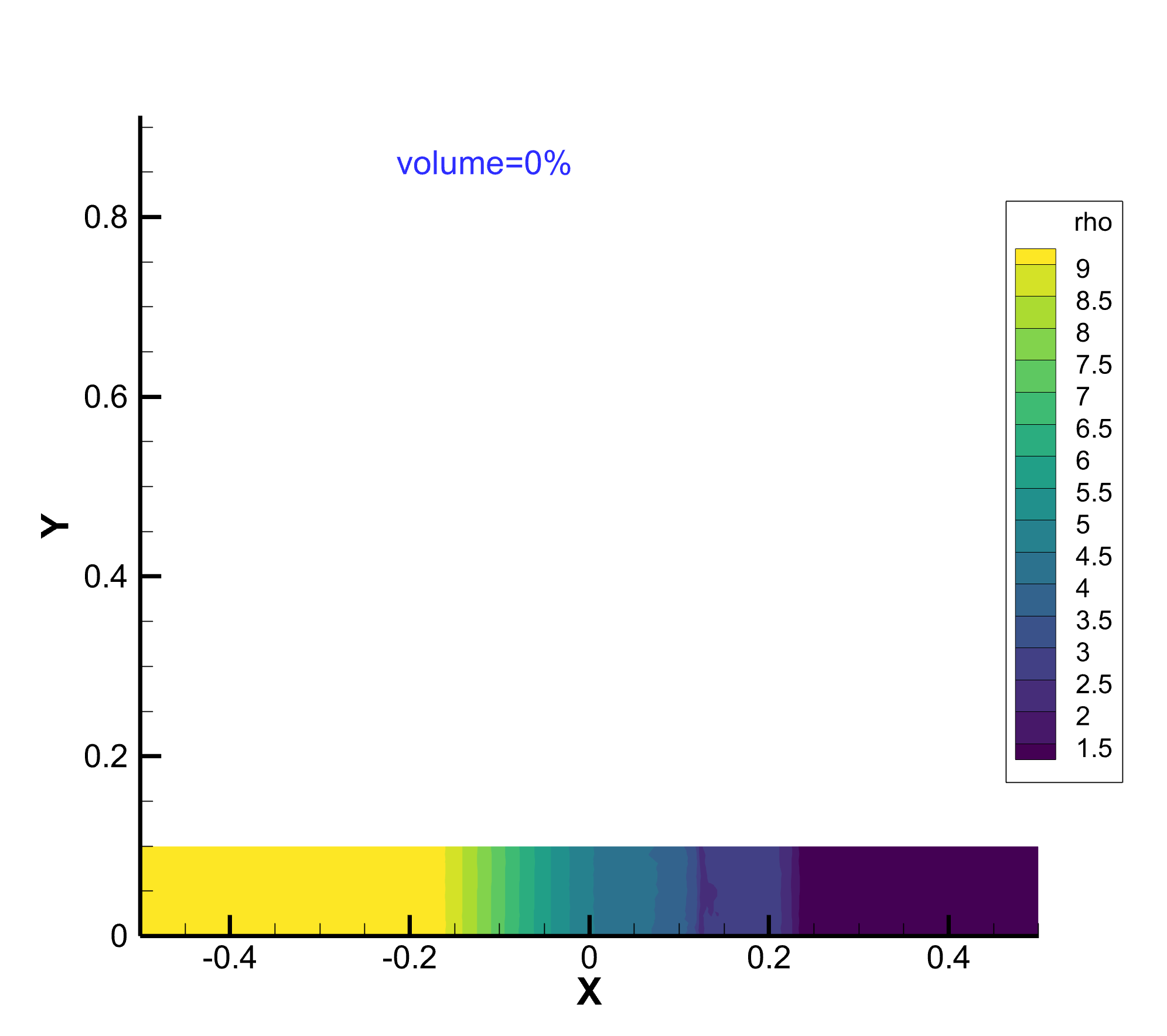}
	\end{minipage}
	\caption{1-D and 2-D Density}   
	\label{fig:density12}    
\end{figure}

\begin{figure}[htbp]
	\begin{minipage}[t]{0.5\textwidth}
		\centering
		\includegraphics[width=\textwidth]{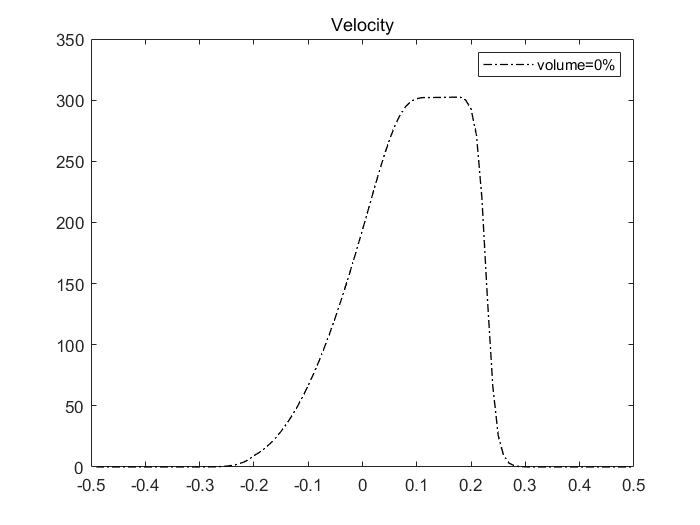}
	\end{minipage}\hfill
	\begin{minipage}[t]{0.45\textwidth}
		\centering
		\includegraphics[width=\textwidth]{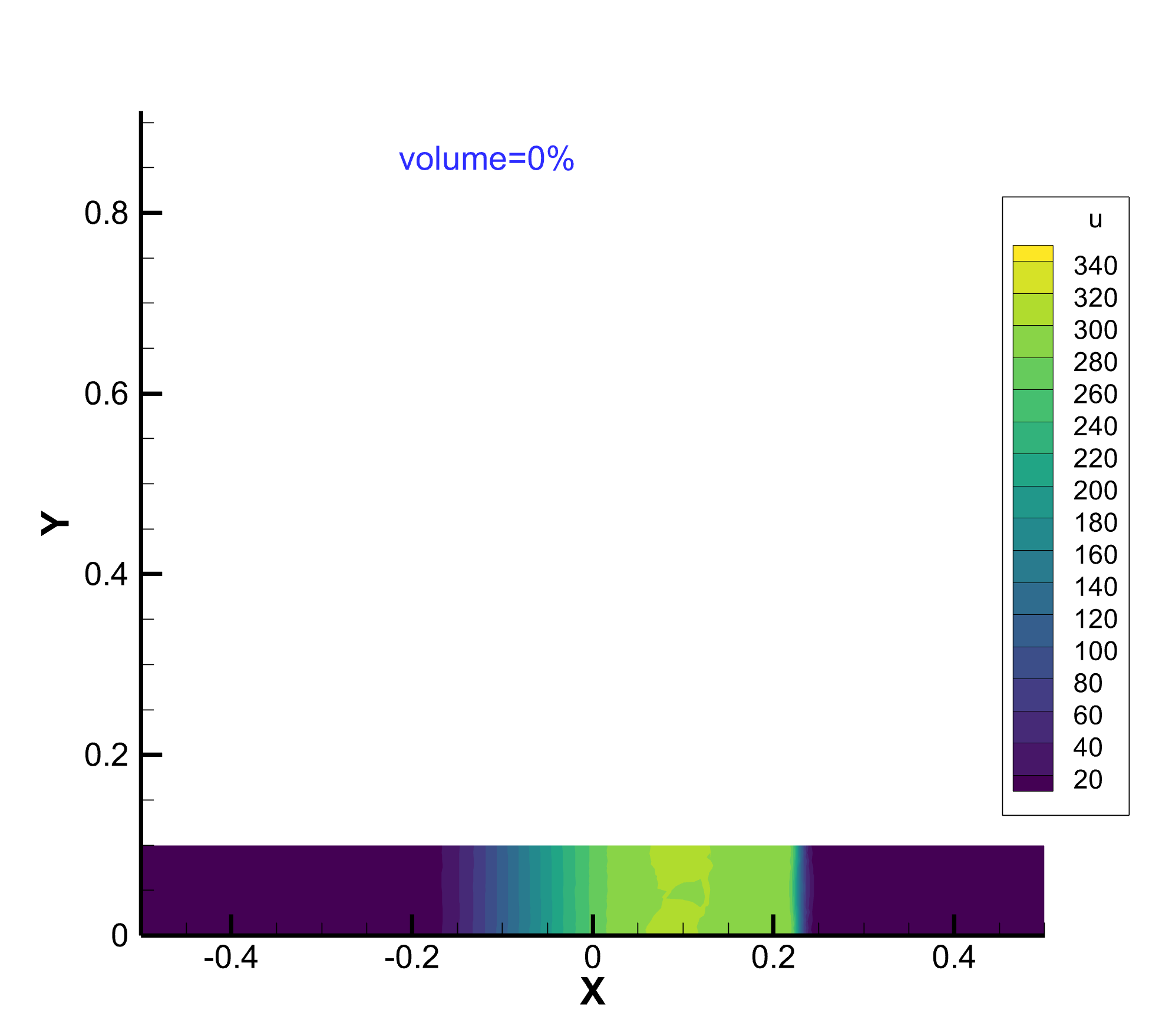}
	\end{minipage}
	\caption{1-D and 2-D Velocity}   
	\label{fig:velocity12}    
\end{figure}

\begin{figure}[htbp]
	\begin{minipage}[t]{0.5\textwidth}
		\centering
		\includegraphics[width=\textwidth]{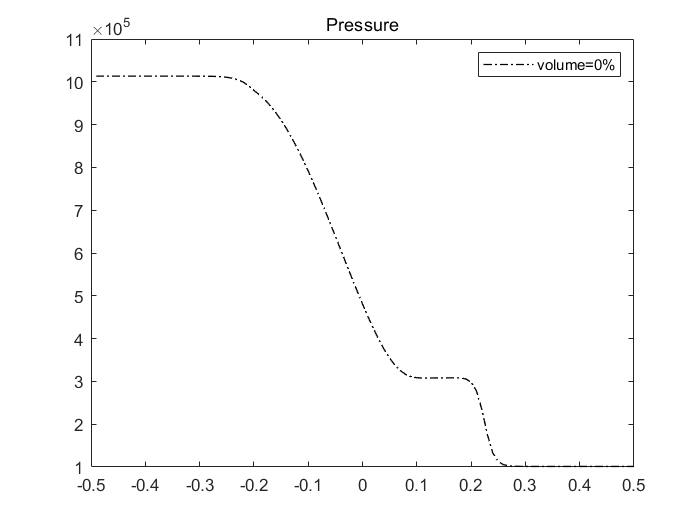}
	\end{minipage}\hfill
	\begin{minipage}[t]{0.45\textwidth}
		\centering
		\includegraphics[width=\textwidth]{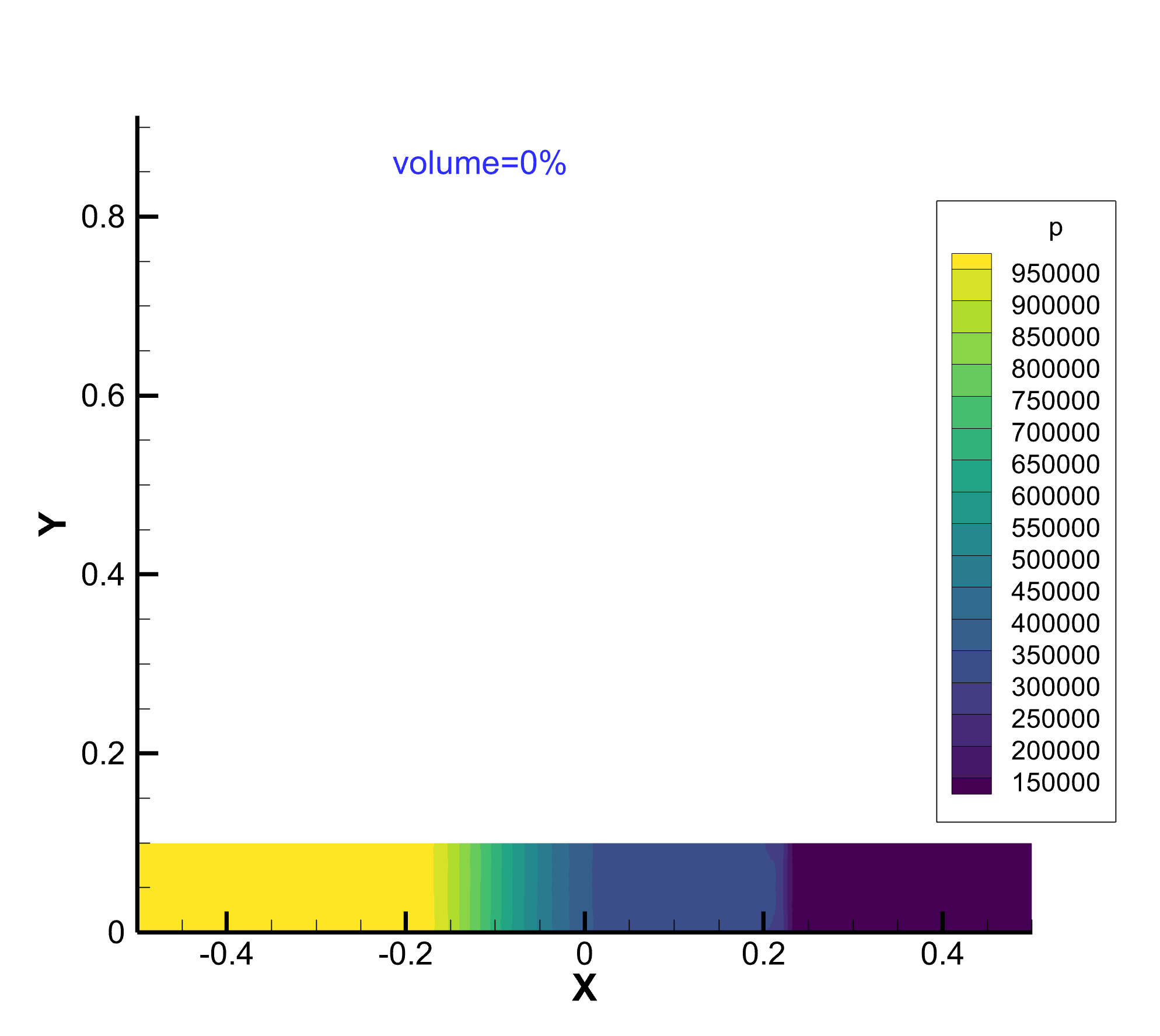}
	\end{minipage}
	\caption{1-D and 2-D Pressure}   
	\label{fig:pressure12}    
\end{figure}

In this tests, we run it in the 2-D region. But for the clearer comparison, we provide both 1-D and 2-D figure. Firstly in figure \ref{fig:density12},\ref{fig:velocity12} and \ref{fig:pressure12} there are the results for the motion of fluid without particles at time step $t=0.0004$. 

\begin{figure}[htbp]
	\begin{minipage}[t]{0.5\textwidth}
		\centering
		\includegraphics[width=\textwidth]{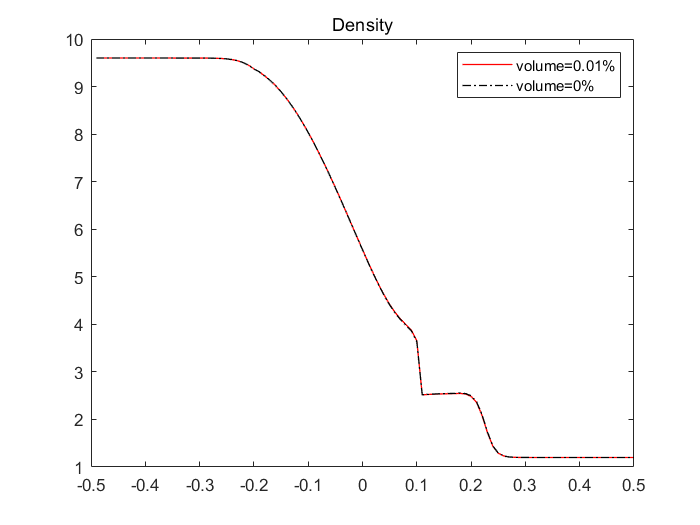}
	\end{minipage}\hfill
	\begin{minipage}[t]{0.45\textwidth}
		\centering
		\includegraphics[width=\textwidth]{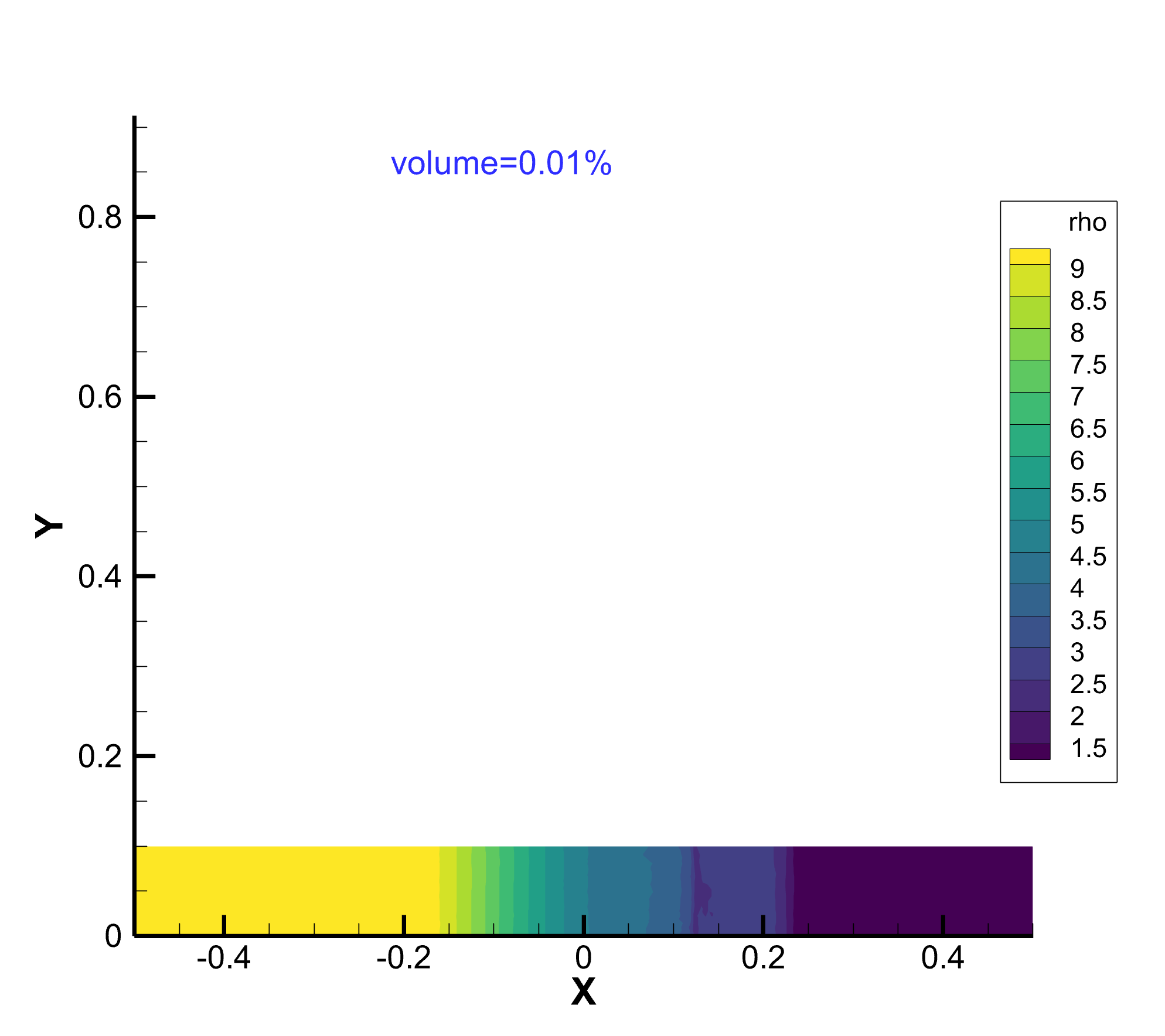}
	\end{minipage}
	\caption{1-D and 2-D Pressure}   
	\label{fig:density_0.01_12}    
\end{figure}

\begin{figure}[htbp]
	\begin{minipage}[t]{0.5\textwidth}
		\centering
		\includegraphics[width=\textwidth]{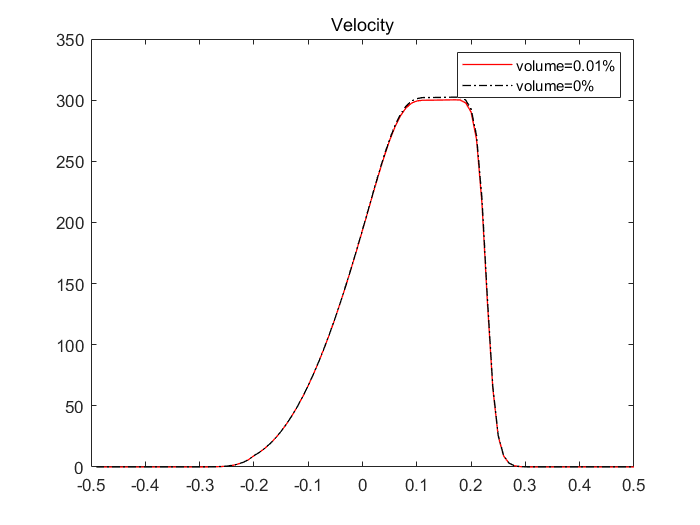}
	\end{minipage}\hfill
	\begin{minipage}[t]{0.45\textwidth}
		\centering
		\includegraphics[width=\textwidth]{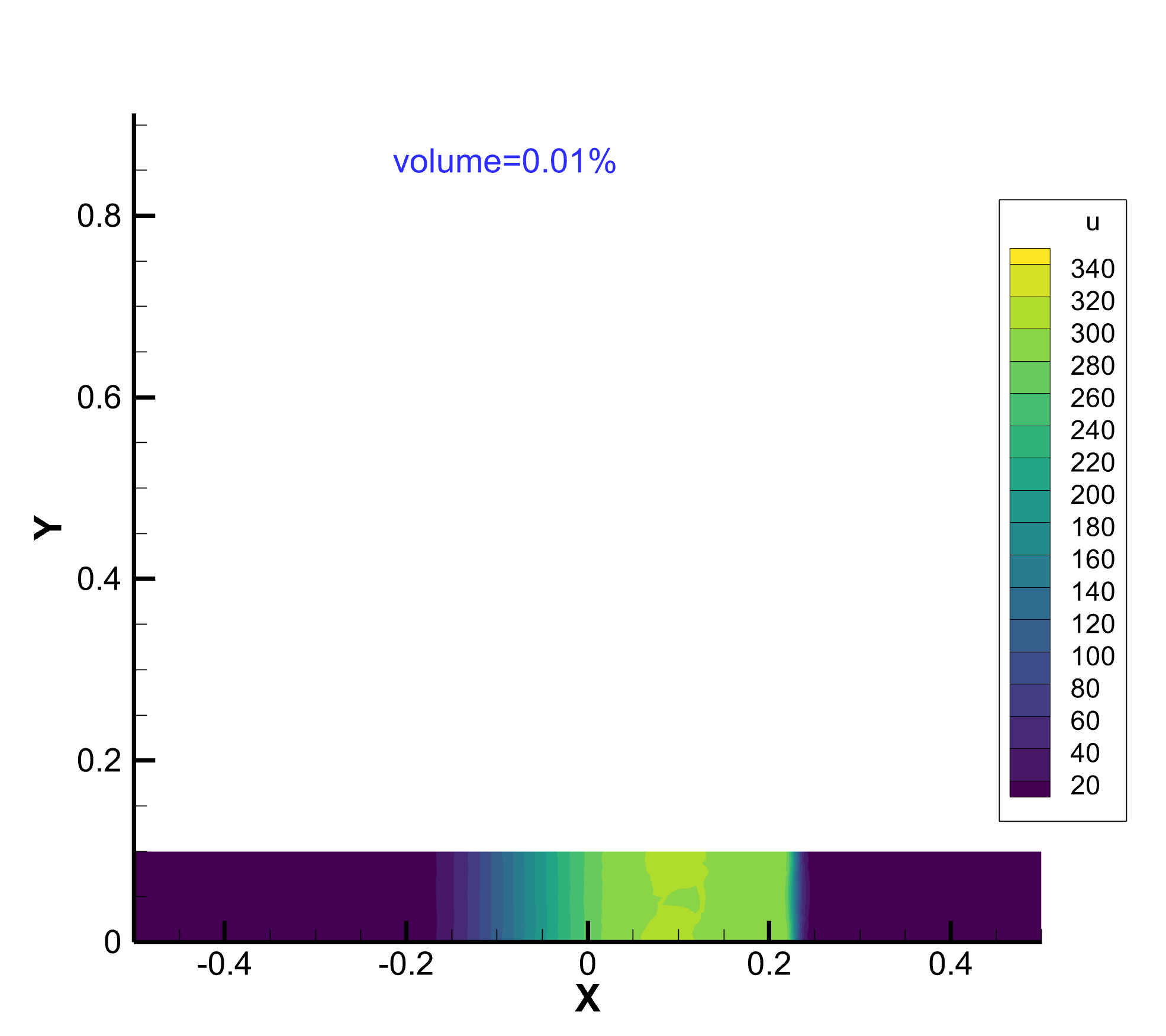}
	\end{minipage}
	\caption{1-D and 2-D Pressure}   
	\label{fig:velocity_0.01_12}    
\end{figure}

\begin{figure}[htbp]
	\begin{minipage}[t]{0.5\textwidth}
		\centering
		\includegraphics[width=\textwidth]{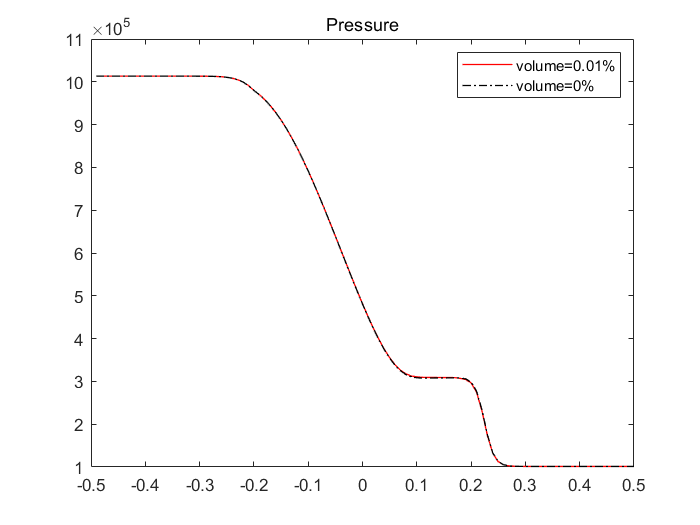}
	\end{minipage}\hfill
	\begin{minipage}[t]{0.45\textwidth}
		\centering
		\includegraphics[width=\textwidth]{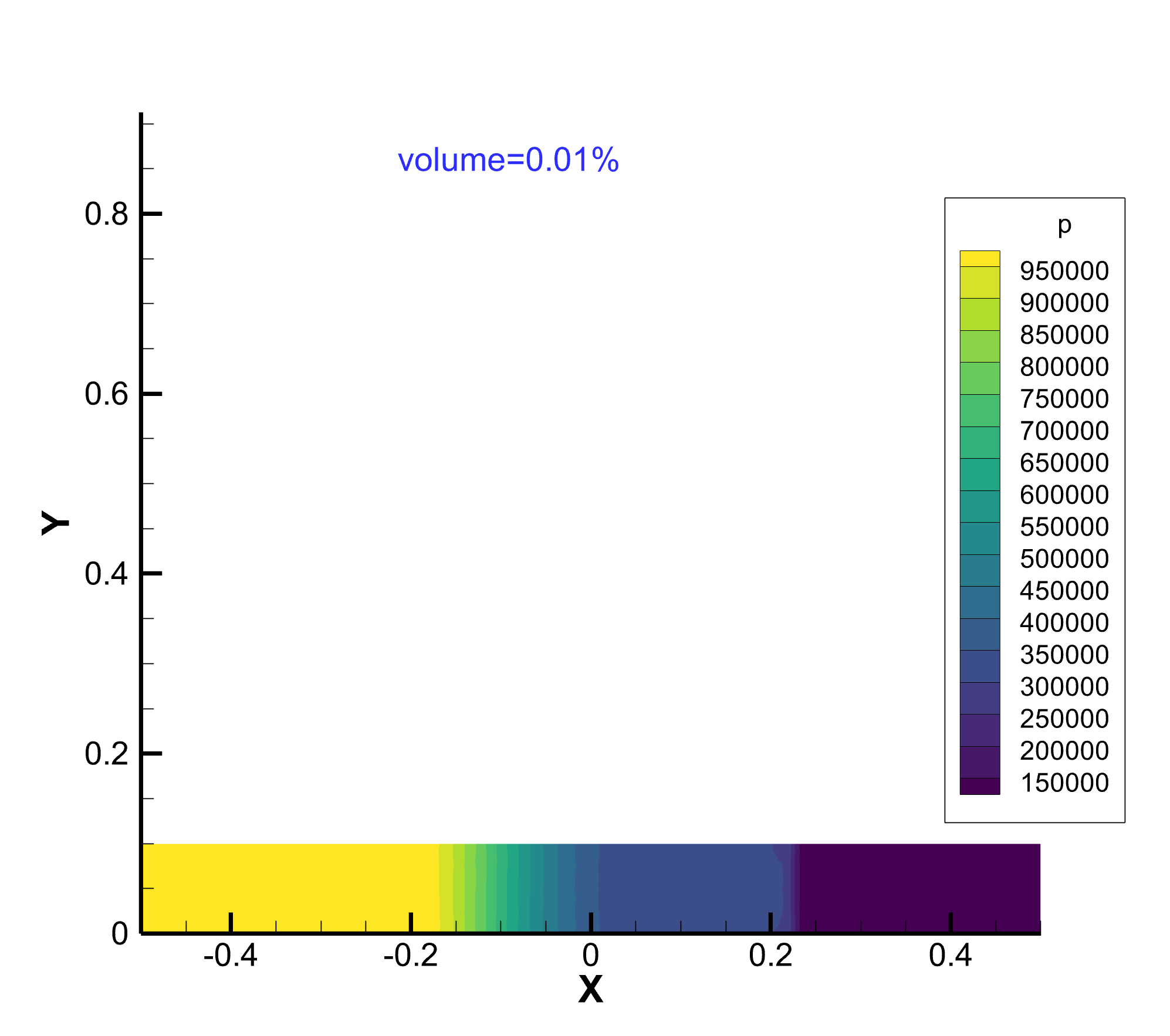}
	\end{minipage}
	\caption{1-D and 2-D Pressure}   
	\label{fig:pressure_0.01_12}    
\end{figure}

Then we change the radius for the particles and give the results in figure \ref{fig:density_0.01_12},\ref{fig:velocity_0.01_12} and \ref{fig:pressure_0.01_12}. From the 1-D results, comparing the black line and red line represented the different the situation of  
particles. We can conclude that when the particle volume fraction is $0.01\%$,  the influence on the fluid is almost negligible. It is difficult to discern especially in density and pressure curves, but a slight effect observed in the velocity curve. This suggests that under the above conditions, particles primarily affect the fluid velocity through force source terms.

Next as we increase the volume percentage of the particles to $4\%$, we can get a significant influence on fluid motion comparing pressure and density which is similar to the results in \cite{tbl}. The results are as follows in figure \ref{fig:density and pressure 4}:

\begin{figure}[htbp]
	\begin{minipage}[t]{0.5\textwidth}
		\centering
		\includegraphics[width=\textwidth]{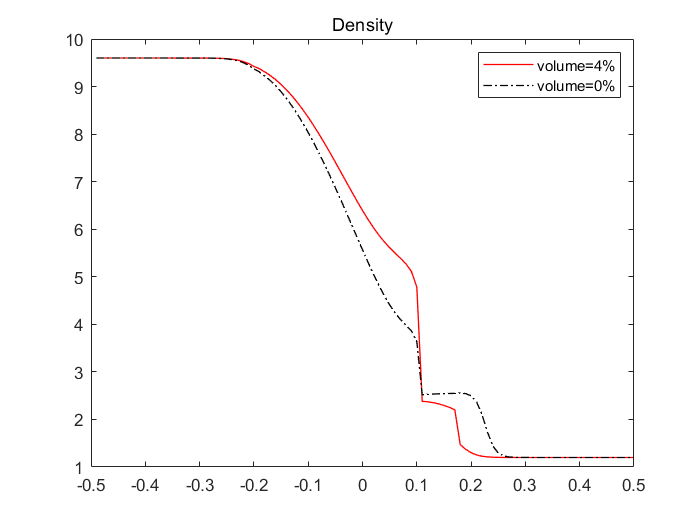}
		Density
	\end{minipage}\hfill
	\begin{minipage}[t]{0.45\textwidth}
		\centering
		\includegraphics[width=\textwidth]{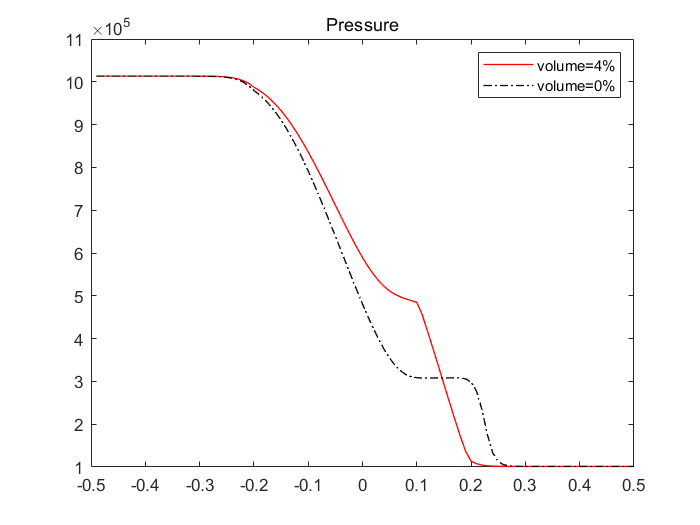}
		Pressure
	\end{minipage}
	\caption{Density and Pressure with $4\%$ particles}   
	\label{fig:density and pressure 4}    
\end{figure}

\begin{remark}
	Although the equations of our model are concerned with the heat exchange between the two phases, we usually ignore it by setting the same temperature in numerical tests. Actually apart from a few examples, energy exchange caused by temperature are not a concern.
\end{remark}

\subsection{Conclusion}
\ 
\newline

Because of the interactions between two phases, the motion of fluid and particles is highly correlated. Under extreme and complex physical conditions, the velocity of particles or fluid could be too high to compute in the Lagrangian form. We design the ALE method to solve the problem which does well in adapting to various situations and tracking particles and fluid.

The above numerical results can to some extent represent the mechanisms of particles influence on fluid. However, it doesn't mean that the same results or influence can be obtained under any conditions. In numerical experiments, through a large amount of numerical tests, we found that the results are highly dependent on the choice fo fluid viscosity coefficient which indicates that the motion trend of fluid and particles is highly sensitive to the fluid viscosity coefficient. We need to analysis the realistic physical background when we do the numerical computations. Additionally, we didn't give the change of internal energy in the above tests because the impact of thermal source term on the fluid is negligible compared with force source term in the case that there is almost no difference between the temperature of fluid and particles. Furthermore, due to the low volume fraction of particles, the fluid and particles have the same pressure at the same position of the region. As the fraction increases, it is necessary to consider the effect of pressure differences on motion in the computational process.



\newpage
\section*{Acknowledgment}

This research was supported by the National Natural Science Foundation of China (grant No.12272059, 11871113, 22341302), the National Key Research and Development Program of China (grant No.2020YFA0713602), and the Key Laboratory of Symbolic Computation and Knowledge Engineering of Ministry of Education of China housed at Jilin University.

\bibliographystyle{siam}

\end{document}